\documentclass[final,3p,times]{elsarticle}

\usepackage{amssymb}
\usepackage{amsmath}
\usepackage{xcolor}

\newif\ifshowrevisions
\showrevisionsfalse
\ifshowrevisions
    \newcommand{\rev}[1]{\textcolor{blue}{#1}}
\else
    \newcommand{\rev}[1]{#1}
\fi

\usepackage{graphicx} 
\usepackage{pgfplots}
\pgfplotsset{compat=1.18}
\usepackage{tikz}
\usetikzlibrary{patterns, matrix, backgrounds, fit} 
\usepackage{nicefrac}
\usepackage{svg}
\usepackage{import}
\usepackage{subcaption}
\usepackage{caption}
\usepackage{float}

\usepackage{amsthm}
\theoremstyle{remark}
\newtheorem*{remark}{Remark}
\usepackage{csquotes}

\usepackage[bookmarks=false]{hyperref} 
\biboptions{sort&compress}   

\definecolor{airforceblue}{rgb}{0.36, 0.54, 0.66}
\definecolor{darkgray}{rgb}{0.66, 0.66, 0.66}
\definecolor{amber}{rgb}{1.0, 0.49, 0.0}
\definecolor{chocolate2267451}{RGB}{226,74,51}
\definecolor{steelblue52138189}{RGB}{52,138,189}
\definecolor{mediumpurple152142213}{RGB}{152,142,213}

\journal{Computer Methods in Applied Mechanics and Engineering
}

\begin{document}

\begin{frontmatter}



\title{Boundary-Level-Constrained Refinement for Suppressing Trimming-Induced High-Frequency Outliers in Explicit Isogeometric Analysis}


\author[inst1,inst2]{Christoph Hollweck}
\author[inst4]{Lukas Leidinger}
\author[inst3]{Stefan Hartmann}
\author[inst1]{Marcus Wagner}
\author[inst2]{Roland Wüchner}

\affiliation[inst1]{organization={OTH Regensburg, Labor Finite-Elemente-Methode},
            addressline={Galgenbergstraße 30}, 
            city={93053 Regensburg},
            country={Germany}}

\affiliation[inst2]{organization={Technische Universität München, Lehrstuhl für Statik und Dynamik},
            addressline={Arcisstr. 21}, 
            city={80333 München},
            country={Germany}}

\affiliation[inst3]{organization={DYNAmore GmbH, an ANSYS Company},
            addressline={Industriestr. 2}, 
            city={70565 Stuttgart},
            country={Germany}}

\affiliation[inst4]{organization={Ansys Italy},%
            addressline={Via G.B. Pergolesi 25}, 
            city={Milano},
            postcode={20124},
            country={Italy}}

\begin{abstract}
In explicit dynamics, the critical time step is governed by the maximum eigenfrequency of the semi-discrete system. 

\rev{In Isogeometric Analysis (IGA), open knot vectors introduce a characteristic spectral boundary effect that can lead to high-frequency outliers and restrict the admissible time step. In the row-sum-lumped setting considered here, extending the computational patch and trimming away the exterior boundary functions mitigates this time-step penalty, but does not necessarily prevent boundary-adjacent basis functions from governing the maximum eigenfrequency.}

\rev{We show that reduced-support basis functions adjacent to the trimming boundary can remain critical spectral contributors even for knot-exact trimming, i.e., in the absence of small cut cells. Hence, small cut cells are not required for trimming-induced time-step penalties. Knot-exact and arbitrary trimming represent different degrees of support reduction within the same underlying mechanism.}

\rev{To control this effect, we propose the \emph{Boundary-Level-Constrained Refinement} (BLCR) strategy, a local refinement constraint for LR- and THB-splines. BLCR constrains the refinement level of basis functions whose support intersects the trimming boundary relative to that of the refined interior. In all configurations investigated in this work, this constraint ensures that the maximum eigenfrequency is governed by untrimmed interior basis functions rather than by refined trimmed functions. Consequently, the trimming-induced high-frequency outliers governing the critical time step are suppressed, yielding a larger admissible time step than globally refined trimmed B-splines with the same interior resolution.}

\rev{BLCR preserves the diagonal row-sum-lumped mass formulation and requires no additional stabilization. The method specifically targets the upper spectral limit governing the critical time step and is not intended to resolve general spectral-accuracy limitations of mass lumping. The underlying mechanism is analyzed using mathematical bounds on the maximum system eigenfrequency and validated in an industrial sheet metal forming simulation with LS-DYNA, where substantially increased stable time steps are obtained while the investigated response quantities remain in close agreement.}
\end{abstract}

\begin{graphicalabstract}

\end{graphicalabstract}

\begin{highlights}
\item \rev{Identifies reduced-support basis functions adjacent to trimming boundaries as critical contributors to the maximum eigenfrequency, even for knot-exact trimming after removal of the exterior open-knot boundary functions}

\item \rev{Explains the spectral penalty caused by reduced support at trimming boundaries using computable eigenfrequency bounds}

\item \rev{Introduces the \emph{Boundary-Level-Constrained Refinement} (BLCR) strategy to suppress the additional time-step penalty caused by refined trimmed basis functions}

\item Achieves critical time step increases of up to 42\% compared to globally refined trimmed patches

\item Preserves the diagonal mass formulation without additional stabilization or mass modification

\item Demonstrates that the BLCR concept is applicable to both LR- and THB-spline discretizations

\item \rev{Demonstrates increased stable time steps in a nonlinear cross-bowl forming simulation with LS-DYNA while maintaining close agreement in the investigated response quantities}

\end{highlights}


\end{frontmatter}


\section{Introduction}
\label{sec:intro}

Isogeometric Analysis (IGA) \cite{hughes_isogeometric_2005} has emerged as a promising framework to bridge the gap between Computer-Aided Design (CAD) and numerical simulation. By utilizing spline-based basis functions directly for analysis, IGA provides superior geometric accuracy and higher-order continuity, which is particularly beneficial for shell structures and complex nonlinear problems \cite{kiendl_isogeometric_2009, bauer_nonlinear_2016}. In industrial practice, geometries are commonly described using Boundary Representation (B-Rep) models. Isogeometric Boundary Representation Analysis (IBRA) \cite{breitenberger_analysis_2015} extends IGA to such models by enabling the analysis of trimmed surfaces and non-conforming couplings, thereby significantly reducing preprocessing effort.

For industrial-scale simulations involving highly nonlinear processes such as crashworthiness or sheet metal forming, explicit time integration is the method of choice. When combined with a diagonal mass matrix, explicit schemes avoid the solution of large linear systems of equations and offer excellent scalability. However, the stable time step is strictly governed by the maximum eigenfrequency of the semi-discrete system. Consequently, high-frequency spectral outliers at the upper end of the discrete spectrum can severely restrict the admissible time step and thereby dominate the computational cost.

\rev{Throughout this work, row-sum mass lumping is treated as a practical efficiency constraint of the targeted industrial explicit dynamics setting. For trimmed discretizations with maximum-continuity splines, Leidinger \cite{leidinger_Diss} showed that row-sum lumping removes the asymptotic dependence of the maximum eigenfrequency on arbitrarily small cut cells, whereas this favorable behavior is not obtained for $C^0$ discretizations. This behavior has subsequently been investigated in greater detail in numerical and theoretical studies of mass lumping for immersed and isogeometric discretizations \cite{Reali_explicit, Bioli2025, Voet_lumping_stabilization}. Complementary approaches based on spectral stabilization, outlier removal, and mass scaling have been developed to further improve the stability and efficiency of explicit immersed and isogeometric discretizations \cite{Voet_outlier, Voet_masscaling, Guarino2025}.}

\rev{The beneficial effect of row-sum mass lumping on the upper end of the spectrum is accompanied by an important stability--accuracy trade-off. Recent studies have shown that, while the largest eigenfrequency remains bounded for sufficiently smooth trimmed spline discretizations, trimming-related spectral pathologies may instead occur in the lower part of the spectrum, where spurious low-frequency modes can adversely affect the transient response if activated \cite{Bioli2025, Voet_lumping_stabilization, Guarino2025}. Radtke et al. \cite{Reali_explicit} observed pronounced accuracy effects in wave-propagation problems, whereas the influence of mass lumping was less pronounced in their nonlinear elastoplastic impact example. These findings indicate a strong dependence on the considered application and nonlinear regime, but do not establish that mass-lumping-related spectral inaccuracies are generally negligible in strongly nonlinear simulations. Their relevance in such applications therefore remains an open question. The present work does not attempt to resolve this general accuracy problem. Instead, it focuses specifically on the upper spectral limit governing the critical explicit time step. Changes observed in the lower part of the spectrum as a consequence of the proposed refinement constraint are, however, explicitly discussed together with their potential implications for solution accuracy.}

\rev{For smooth $C^{p-1}$ spline discretizations, an additional limitation of the critical time step is associated with the boundaries of clamped, or open, knot vectors. Adam et al. \cite{ADAM2015581} showed that the corresponding boundary functions can govern the stable time step and investigated two remedies. Increasing the size of the boundary elements can shift the time-step limitation towards interior functions, but the required coarsening becomes substantial and increasingly degree dependent and may compromise approximation quality when boundary effects are relevant. As an alternative, Adam et al. proposed unclamping the knot vector, thereby removing the particularly restrictive exterior functions of the clamped representation while preserving the underlying geometry through a corresponding reparameterization.}

\rev{In trimmed B-Rep discretizations, a related approach is to choose the background spline patch sufficiently large such that the exterior open-knot boundary spans lie outside the physical domain and are removed by trimming \cite{leidinger_Diss, mesmer_efficient_2022, hollweck_LR_THB_2026}. Me{\ss}mer et al. \cite{mesmer_efficient_2022}, for example, employ this construction to obtain feasible critical time steps for trimmed high-continuity B-spline discretizations and identify non-open knot vectors as an alternative. In both approaches, the particularly restrictive exterior open-knot boundary functions are absent from the active representation of the physical domain.}

\rev{Removing these exterior functions, however, does not eliminate all trimming-related contributions to the maximum eigenfrequency. Knot-exact and arbitrary trimming should not be regarded as fundamentally different mechanisms in this respect. In both cases, trimming restricts the active portion of the support of adjacent basis functions and modifies their stiffness and row-sum-lumped mass contributions non-proportionally. The resulting stiffness-to-mass ratio depends on the particular trimming configuration and on which portion of the basis function support remains active. Consequently, an arbitrarily small cut fraction is not required for a reduced-support basis function to become relevant for the maximum eigenfrequency. Knot-exact trimming therefore provides a particularly useful configuration for isolating this support-reduction mechanism from effects associated exclusively with vanishingly small cut fractions.}

\rev{Local refinement introduces a second, conceptually distinct reduction of the relevant length scale. While trimming restricts the active portion of a given basis-function support, refinement generates basis functions with a smaller characteristic support size. If such refined basis functions are additionally intersected by the trimming boundary, both effects act simultaneously: the nominal support is reduced by refinement and its active portion is further restricted by trimming. How this interaction affects the upper spectral limit is the central question addressed in the present work.}

In large-scale explicit simulations, local refinement is essential to efficiently resolve localized phenomena such as plastic deformation, contact-induced stress concentrations, or fracture. Established local refinement techniques in IGA include Locally Refined (LR) splines \cite{dokken_polynomial_2013} and Truncated Hierarchical B-splines (THB) \cite{giannelli_thb-splines_2012}.

\rev{While the boundary effect of open knot vectors and corresponding remedies are well established, it remains insufficiently understood how the reduced active support caused by trimming interacts with the smaller characteristic support introduced by local refinement. In particular, refinement of basis functions whose active support is already restricted by trimming may cause these functions to become the governing contributors to the maximum eigenfrequency. The present work systematically investigates this interaction and identifies reduced-support basis functions adjacent to the trimming boundary as the relevant contributors to this remaining time-step limitation.}

\rev{Based on this observation, we introduce the \emph{Boundary-Level-Constrained Refinement} (BLCR) strategy, which constrains the refinement level of basis functions whose support intersects the trimming boundary relative to the refined interior. BLCR specifically targets the upper spectral limit governing the critical time step and introduces an application-dependent trade-off between boundary resolution and time-step efficiency. It is therefore intended primarily for applications in which local refinement is required in the interior while the trimmed boundary is not a region of primary approximation interest. BLCR should be understood as a targeted strategy for controlling the critical time step rather than as a general remedy for spectral-accuracy limitations associated with mass lumping.}

\subsection*{Aim and Outline of the Paper}

\rev{Building on our previous work \cite{hollweck_LR_THB_2026}, the present study investigates the effect of trimming and local refinement on the upper spectral limit of row-sum-lumped LR- and THB-spline discretizations and develops BLCR as a corresponding refinement strategy.}

\rev{The main contributions of this paper are:}
\begin{itemize}
    \item \rev{Identification of reduced-support trimmed basis functions as the remaining critical contributors to the maximum eigenfrequency.}
    
    \item \rev{A bound-based explanation of their contribution to the upper spectral limit for knot-exact and arbitrary trimming configurations.}
    
    \item \rev{The formulation and validation of BLCR for LR- and THB-spline discretizations, yielding larger admissible time steps than globally refined trimmed patches at the same interior resolution.}
\end{itemize}

\rev{The remainder of the paper is organized as follows. Section~\ref{sec:sec2} introduces the spline technologies, trimming concepts, stability criteria, and numerical models. Section~\ref{sec:sec3} presents the spectral analyses and the BLCR strategy. Section~\ref{sec:sec4} evaluates BLCR in an industrial nonlinear explicit benchmark and discusses its applicability and limitations. Section~\ref{sec:sec5} summarizes the main findings and outlines directions for future work.}

\section{Theoretical Framework}
\label{sec:sec2}

This section provides a concise introduction to B-spline formulations, which serve as the foundation for LR- and THB-splines. While this study focuses on B-splines, a comprehensive overview of Non-Uniform Rational B-splines (NURBS) is provided in \cite{piegl_nurbs_1997}. Both LR- and THB-splines can be extended to a rational basis, as detailed in \cite{Karsten_2019, Grendas_2023}. Furthermore, we introduce the concept of trimming and the definition of Boundary Representation (B-Rep) models as used in standard CAD environments. The section concludes with the theoretical stability limits for explicit time integration and the numerical models utilized for the subsequent analysis.
\\

\subsection{B-splines}

\subsubsection*{Univariate B-splines}

Let $\tilde{\Omega}\subset\mathbb{R}$ denote a one-dimensional
parametric domain with coordinate $\xi$.
Let $n$ denote the number of basis functions.
Given a knot vector
\[
\Xi=\{\xi_1,\xi_2,\dots,\xi_{n+p+1}\},
\]
B-spline basis functions $N_{i,p}(\xi)$ of degree $p$
are defined via the Cox--de Boor recursion \cite{DEBOOR197250}.
For $\xi_i\le\xi<\xi_{i+1}$,
\[
N_{i,0}(\xi)=1,
\]
and zero otherwise. Higher degrees are defined recursively as
\[
N_{i,p}(\xi)
=
\frac{\xi-\xi_i}{\xi_{i+p}-\xi_i}N_{i,p-1}(\xi)
+
\frac{\xi_{i+p+1}-\xi}{\xi_{i+p+1}-\xi_{i+1}}N_{i+1,p-1}(\xi).
\]
The knot spans
\begin{equation}
    e_i = [\xi_i, \xi_{i+1}]
\end{equation}
with $\xi_i < \xi_{i+1}$ define a partition of $\tilde{\Omega}$ into elements of non-zero length.
B-splines are piecewise polynomials of degree $p$ with $C^{p-m}$ 
continuity at knots of multiplicity $m$, and reduce to standard 
polynomials of degree $p$ on each element. They form a partition of unity,
\[
\sum_i N_{i,p}(\xi)=1,
\]
and have compact support
\[
\operatorname{supp}(N_{i,p})=[\xi_i,\xi_{i+p+1}].
\]

\subsubsection*{Tensor-product B-splines}

In two dimensions, the parametric domain is defined as
\[
\tilde{\Omega}
=
[\xi_{\min},\xi_{\max}]
\times
[\eta_{\min},\eta_{\max}]
\subset\mathbb{R}^2,
\]
with parametric coordinates $\boldsymbol{\xi}=(\xi,\eta)$.
Given knot vectors in both directions,
tensor-product basis functions are defined as
\[
N_A(\xi,\eta)
=
N_{i,p}(\xi)\,M_{j,q}(\eta),
\]
\rev{where $A$ denotes the global basis-function index associated with the tensor-product index pair $(i,j)$.}
The parametric mesh is given by the tensor-product partition
\begin{equation}
    e_{ij} = [\xi_i, \xi_{i+1}] \times [\eta_j, \eta_{j+1}] \subset \tilde{\Omega},
\end{equation}
where $\xi_i < \xi_{i+1}$ and $\eta_j < \eta_{j+1}$. Each element $e_{ij}$ forms a rectangle in the parametric domain on which the bivariate basis functions are polynomials of degree $(p, q)$.
A spline surface is defined by its control points 
$\mathbf{P}_A \in \mathbb{R}^3$ as
\[
\mathbf{S}(\xi,\eta)
=
\sum_A N_A(\xi,\eta)\,\mathbf{P}_A,
\]
which induces the geometric mapping
\[
\mathbf{S}:\tilde{\Omega}\rightarrow\Omega\subset\mathbb{R}^3.
\]
Hence, spline basis functions are defined on the parametric domain
$\tilde{\Omega}$, while the physical geometry $\Omega$
is obtained as its image under $\mathbf{S}$.
Throughout this work, we assume $p=q$.

B-splines can be extended to Non-Uniform Rational B-Splines (NURBS),
which constitute the industry standard in CAD due to their ability
to represent conic sections exactly \cite{piegl_nurbs_1997}.
However, many engineering applications can be accurately approximated
using polynomial B-splines.
In the present study, standard B-splines are employed in order
to reduce computational complexity while maintaining sufficient
geometric accuracy for the problems under consideration.

\subsection{LR-splines}
\label{sec:lrsplines}

LR-splines are based on the concept of local knot insertion.
We first introduce the univariate refinement procedure and then extend it
to the bivariate case. For comprehensive treatments we refer to
\cite{dokken_polynomial_2013, johannessen_isogeometric_2014, hollweck_LR_THB_2026}.

\subsubsection*{Univariate Refinement}

We consider a one-dimensional parametric domain
\[
\tilde{\Omega} = [\xi_1,\xi_{n+p+1}],
\]
equipped with a knot vector
\[
\Xi = \{\xi_1,\dots,\xi_{n+p+1}\}.
\]
Again, the knot vector induces a partition of $\tilde{\Omega}$ into
non-degenerate knot spans
\[
e_i = [\xi_i,\xi_{i+1}],
\qquad
\xi_i < \xi_{i+1}.
\]
The collection of all such intervals defines the initial mesh
$\mathcal{M}^0$.

Consider a basis function $N_{\Xi_i}$ associated with the local knot vector
$\Xi_i = \{\xi_i,\dots,\xi_{i+p+1}\}$.
Upon inserting a knot $\bar{\xi}$ with
$\xi_i \le \bar{\xi} \le \xi_{i+p+1}$,
the parent basis function is represented as a weighted sum of two children,
\begin{equation}
N_{\Xi_i}(\xi)
=
\alpha_1 N_{\Xi_1^*}(\xi)
+
\alpha_2 N_{\Xi_2^*}(\xi),
\label{eq:funsplit}
\end{equation}
where
\begin{equation}
\alpha_1 =
\begin{cases}
\rev{\dfrac{\bar{\xi}-\xi_i}{\xi_{i+p}-\xi_i}},
& \rev{\xi_i \le \bar{\xi} < \xi_{i+p}}, \\[6pt]
\rev{1},
& \rev{\xi_{i+p} \le \bar{\xi} \le \xi_{i+p+1}},
\end{cases}
\qquad
\alpha_2 =
\begin{cases}
\rev{1},
& \rev{\xi_i \le \bar{\xi} \le \xi_{i+1}}, \\[4pt]
\rev{\dfrac{\xi_{i+p+1}-\bar{\xi}}{\xi_{i+p+1}-\xi_{i+1}}},
& \rev{\xi_{i+1} < \bar{\xi} \le \xi_{i+p+1}}.
\end{cases}
\label{eq:alphas}
\end{equation}
Insertion of $\bar{\xi}$ locally subdivides the affected knot spans,
thereby refining the mesh $\mathcal{M}$ in the parametric domain.
A sequence of knot insertions produces a nested mesh hierarchy
\[
\mathcal{M}^k \subset \mathcal{M}^{k+1}.
\]

\subsubsection*{Bivariate Refinement}

In two dimensions, the parametric domain $\tilde{\Omega} \subset \mathbb{R}^2$ is defined by two knot vectors $\Xi$ and $\mathcal{H}$. The tensor product of their knot spans defines the rectangular elements
\begin{equation}
    e_{ij} = [\xi_i, \xi_{i+1}] \times [\eta_j, \eta_{j+1}],
\end{equation}
where $\xi_i < \xi_{i+1}$ and $\eta_j < \eta_{j+1}$. The parametric mesh $\mathcal{M}$ is the collection of all such non-degenerate rectangles.
A bivariate basis function is defined by two local knot vectors
$\Xi_i$ and $\mathcal{H}_j$ and carries a weight $\gamma_A$,
\[
N_{\mathbf{\Xi}_A}^{\gamma_A}(\xi,\eta)
=
\gamma_A N_{\Xi_i}(\xi) N_{\mathcal{H}_j}(\eta).
\]
Initially, each function carries a unit weight $\gamma_A = 1$.
Refinement follows the univariate splitting logic,
\begin{equation}
\begin{aligned}
N_{\mathbf{\Xi}_A}^{\gamma_A}(\xi,\eta)
&=
\gamma_A
\left(
\alpha_1 N_{\Xi_1^*}(\xi)
+
\alpha_2 N_{\Xi_2^*}(\xi)
\right)
N_{\mathcal{H}_j}(\eta)
\\
&=
N_{\mathbf{\Xi}_1^*}^{\gamma_1}(\xi,\eta)
+
N_{\mathbf{\Xi}_2^*}^{\gamma_2}(\xi,\eta),
\end{aligned}
\label{eq:bivariate}
\end{equation}
with $\gamma_1=\alpha_1\gamma_A$ and
$\gamma_2=\alpha_2\gamma_A$.
Local refinement is performed by inserting meshlines in
$\tilde{\Omega}$.
A meshline is an axis-aligned segment defined by fixing one
parametric coordinate and specifying a finite span in the other.
Its endpoints must coincide with existing meshlines.
All elements intersected by the meshline are subdivided
into two child elements by splitting them along the meshline.
The resulting elements remain axis-aligned rectangles.
This produces a locally refined partition of $\tilde{\Omega}$,
leading to a nested sequence of meshes
\[
\mathcal{M}^k \subset \mathcal{M}^{k+1}.
\]
For notational consistency with the subsequent sections,
we denote
\[
T_A := N_{\mathbf{\Xi}_A}^{\gamma_A},
\]
where the resulting set of LR basis functions is denoted by
\[
\mathcal{T}
=
\{ T_A \}.
\]

\subsubsection*{Structured Refinement Strategy}

To ensure comparability with THB-splines, we employ a structured refinement strategy.
An initial tensor-product mesh $\mathcal{M}^0$ is given.
An element $e$ is said to be of level $l$ if it results from
$l$ successive uniform bisections of an initial element.
Each bisection consists of inserting one horizontal and one
vertical meshline through the midpoint of the element,
thereby subdividing it into four child elements.
The child elements are assigned level $l+1$.

Refinement is triggered by marking a basis function $T_A^l$ of level $l$.
Structured refinement inserts meshlines such that all elements in the support of $T_A^l$
are split in each parametric direction such that all elements are at least of level $l+1$. Consequently, all affected elements become level $l+1$ elements. The newly generated basis functions are also assigned level $l+1$. If a basis function has support on elements of different levels,
its level is defined as
\[
\ell(T_A) := \min_{e \subset \operatorname{supp}(T_A)} \ell(e),
\]
that is, the minimum level of all elements contained in its support. This refinement strategy preserves a constant aspect ratio of all elements.

\begin{figure}[h]
    \centering
    \begin{minipage}{0.4\textwidth}
        \centering
        \includegraphics[width=\linewidth]{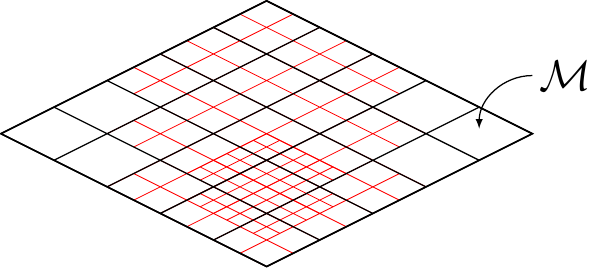}
    \end{minipage}
    \hspace{0.05\textwidth}
    \begin{minipage}{0.45\textwidth}
        \centering
        \includegraphics[width=0.6\linewidth]{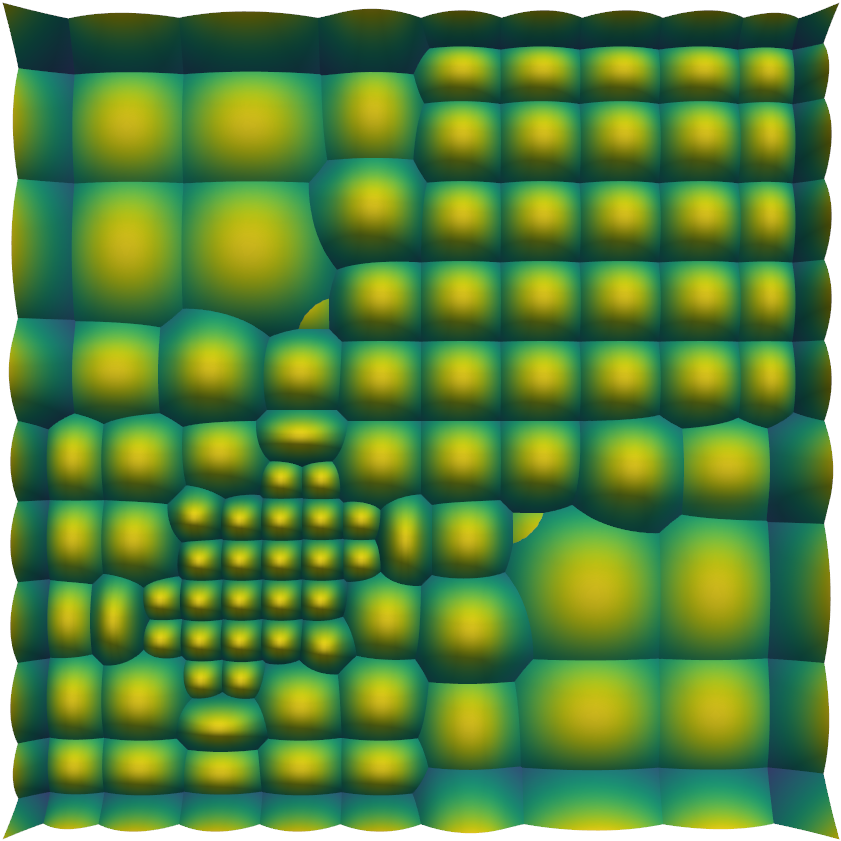}\\[0.5em]
        \includegraphics[width=0.8\linewidth]{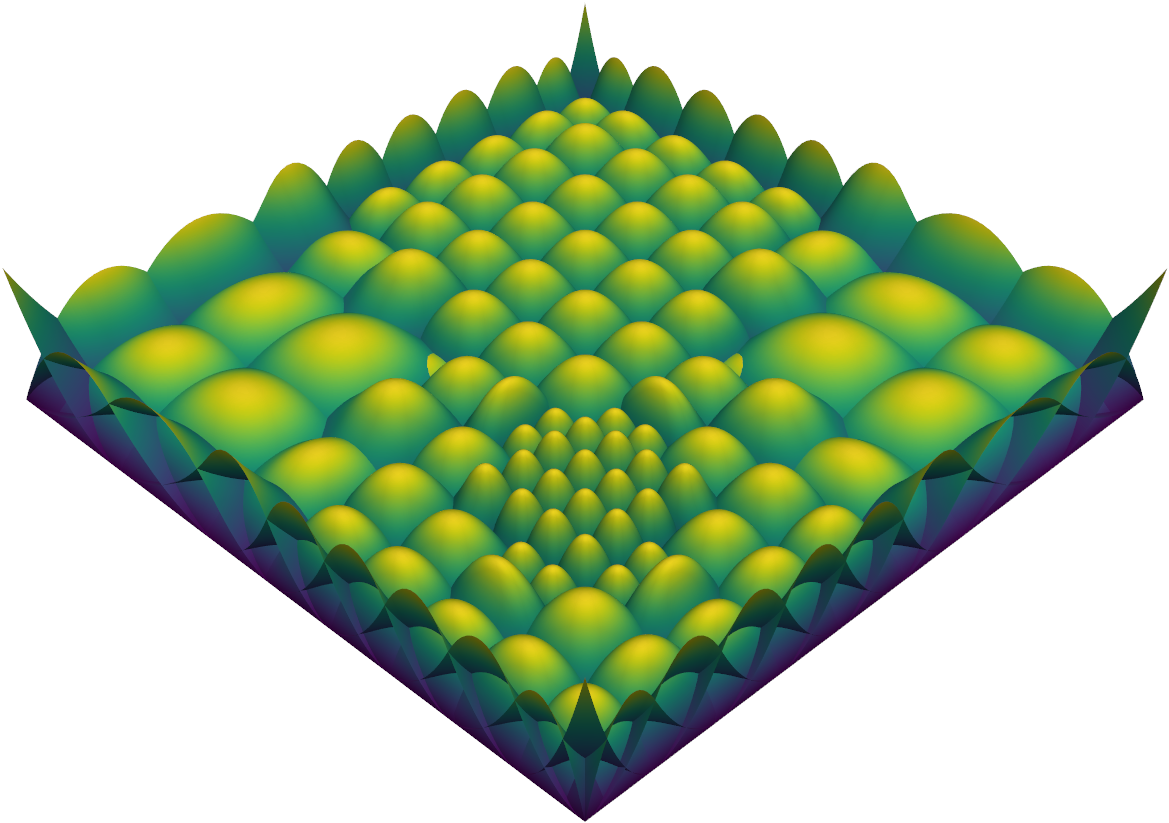}
    \end{minipage}
    \caption{LR refinement: resulting Mesh $\mathcal{M}$ (left) and corresponding bivariate $C^1$ basis functions of degree $p=2$ (right).}
    \label{fig:LR_combined}
\end{figure}

Fig.~\ref{fig:LR_combined} displays a locally refined mesh with the inserted meshlines highlighted in red. Additionally, the final set of bivariate basis functions is shown, demonstrating the localized nature of the refinement. 

\subsection{THB-splines}
\label{sec:thbsplines}

Unlike LR-splines, Truncated Hierarchical B-splines (THB-splines) utilize a hierarchical structure of nested spaces $\mathcal{V}^0 \subset \mathcal{V}^1 \subset \dots \subset \mathcal{V}^L$ to achieve local refinement. Each level $l$ is associated with a mesh $\mathcal{M}^l$ and a set of basis functions $\mathcal{N}^l$. Refinement regions $\tilde{\Omega^{l}}$ are defined such that $\tilde{\Omega}^{l+1} \subseteq \tilde{\Omega}^l$, where active elements are selected across different levels to form the hierarchical mesh $\hat{\mathcal{M}}$:
\begin{equation}
\mathcal{\hat{M}} := \left\{ e \in \mathcal{M}^l : e \subseteq \tilde{\Omega}^l \land e \not\subseteq \tilde{\Omega}^{l+1}, \; l = 0, \dots, L - 1 \right\}.
\label{eq:meshthb}
\end{equation}
The THB-basis $\mathcal{T}$ is defined recursively, incorporating a truncation operator to preserve the partition of unity:
\begin{equation}
\begin{cases}
\mathcal{T}^0 := \mathcal{N}^0, \\

\mathcal{T}^{l+1} := \left\{ \text{trunc}(T_i^l): T_i^l \in \mathcal{T}^l \land \text{supp}(T_i^l) \nsubseteq \tilde{\Omega}^{l+1} \right\} \cup \\

\quad \quad \quad \, \, \left\{ \text{child}(T_i^l)_j \in \mathcal{N}^{l+1} : \text{supp}(T_i^l) \subseteq \tilde{\Omega}^{l+1} \right\} \quad l = 0, \dots, L - 1, \\
\mathcal{T} := \mathcal{T}^L
\end{cases}
\label{eq:recursive}
\end{equation}
The truncation operator modifies coarse basis functions by removing the contribution of child functions that are already active on finer levels:
\begin{equation}
\text{trunc}(T_i^l) = \sum_{\substack{j : \text{child}(T_i^l)_j \in \mathcal{N}^{l+1} \\ \text{supp}(\text{child}(T_i^l)_j) \nsubseteq \tilde{\Omega}^{l+1}}} c_{ij} \cdot \text{child}(T_i^l)_j,
\label{eq:trunc}
\end{equation}
where $c_{ij}$ are the coefficients from the two-scale relation. For a detailed 1D illustration of this truncation procedure and the resulting basis shapes, we refer the reader to \cite{hollweck_LR_THB_2026} and \cite{giannelli_thb-splines_2012}.

\begin{figure}[h]
  \centering
  \parbox{0.3\textwidth}{
    \includegraphics[width=\linewidth]{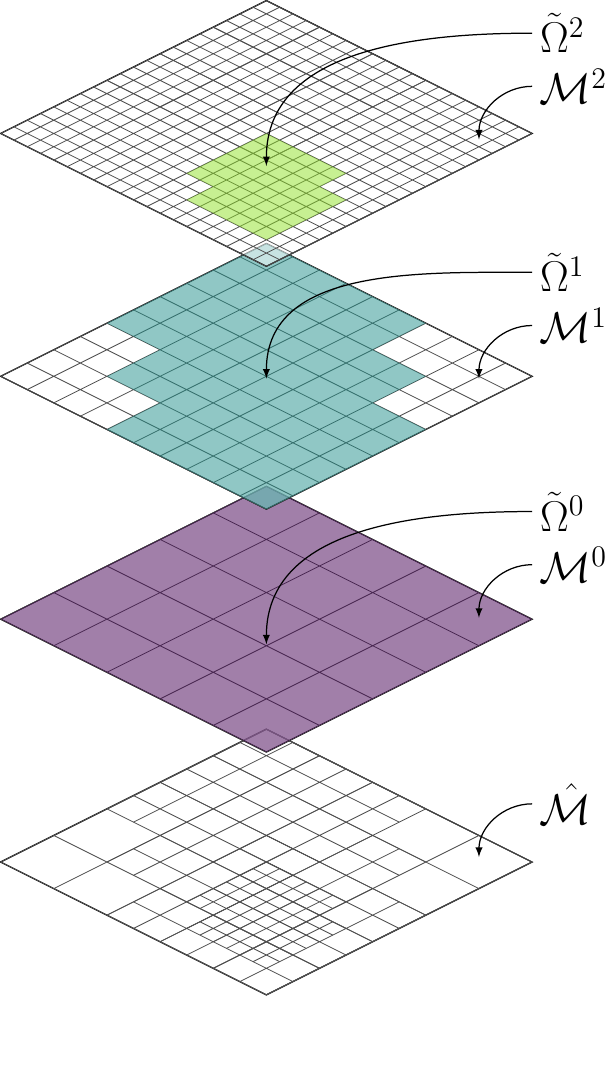}
  }
  \hspace{0.05\textwidth}
  \parbox{0.45\textwidth}{
    \centering
    \includegraphics[width=0.6\linewidth]{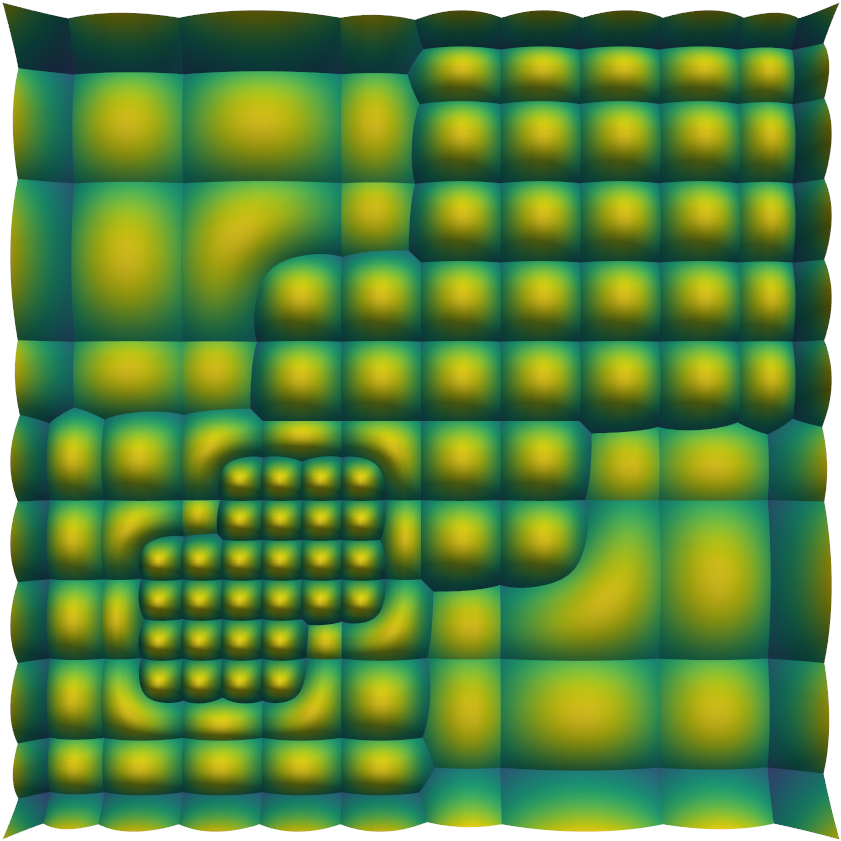}
    \\[0.5em] 
    \includegraphics[width=0.8\linewidth]{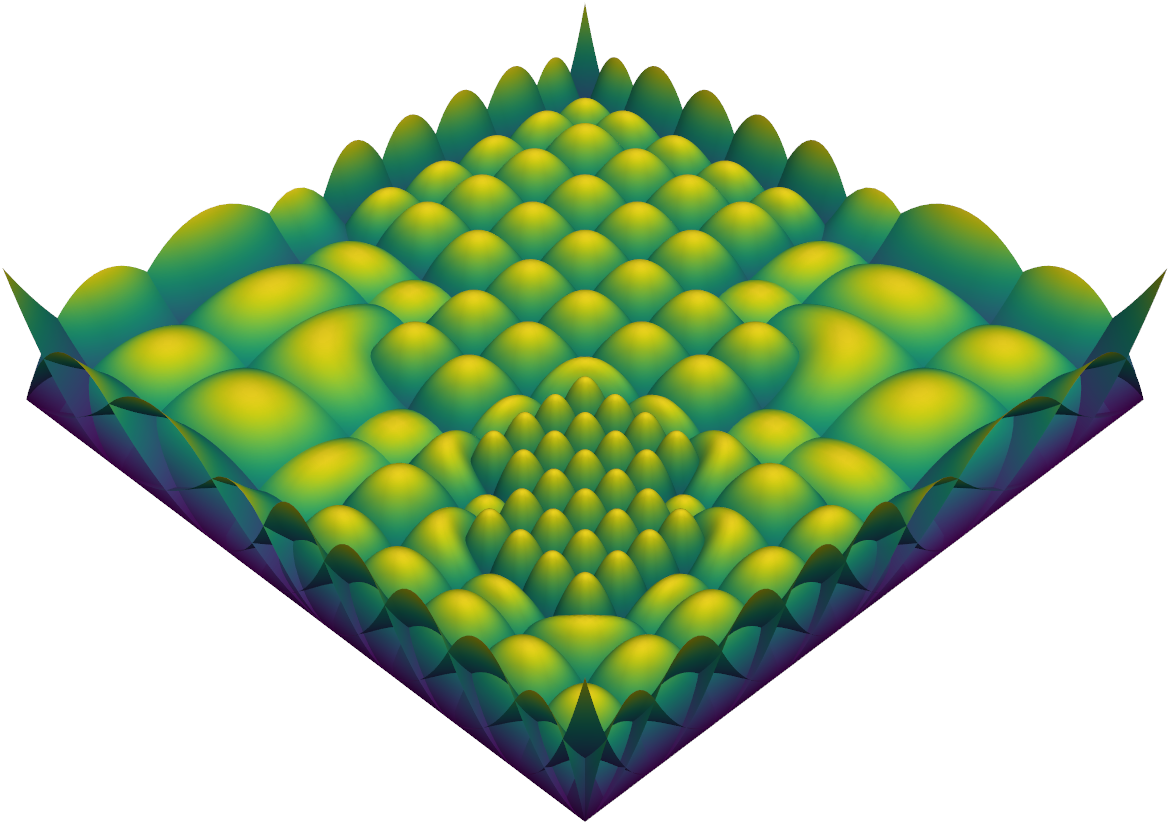}
  }
  \caption{THB-spline refinement: Hierarchical mesh $\hat{\mathcal{M}}$ (left) and the resulting bivariate $C^1$ basis of degree $p=2$ (right).}
  \label{fig:THB_combined}
\end{figure}

Refinement in 2D follows the same recursive principles, though truncation may result in non-convex supports. To facilitate a direct comparison between LR- and THB-splines, the hierarchical mesh shown in Fig.~\ref{fig:THB_combined} is identical to the one utilized for LR-splines in Fig.~\ref{fig:LR_combined}. The illustration highlights the active basis functions and the resulting mesh $\hat{\mathcal{M}}$. In contrast to LR-splines, THB-splines require an explicit update of control points for newly activated functions during adaptive refinement \cite{GIANNELLI2016}. It should also be noted that basis functions can undergo multiple truncations, particularly in regions with sharp transitions between refinement levels.

In the original THB-spline framework, hierarchical refinement is controlled by the support of the basis functions, that is, activation and truncation are determined by geometric overlap of supports \cite{giannelli_thb-splines_2012, giannelli_strongly_2014, Buffa2022_AdaptiveIGA}. The approach employed here follows a different principle. Basis functions are managed through explicit parent--child relationships rather than through support-based criteria. This construction is aligned with the concept presented in \cite{BuffaGarau2016RefinableSpaces}, which builds on earlier strategies developed in \cite{Grinspun2002CHARMS, Krysl2003NaturalHierarchical}. Because the hierarchical dependencies are known in advance, the refinement logic and the associated data structures are straightforward. Each child function can be assigned directly to its generating parent, which avoids additional geometric checks during activation or truncation. From a functional analytic point of view, the resulting discrete space preserves the approximation properties of the classical THB formulation, while typically involving fewer active degrees of freedom. Further algorithmic aspects are documented in \cite{GARAU201858}.

\subsection{Trimmed Surfaces}
\label{sec:trimmedsurfaces}

In industrial CAD environments, geometries are typically represented as Boundary Representation (B-Rep) models, where objects are defined by their outer skins. A B-spline surface $\mathbf{S} : \tilde{\Omega} \rightarrow \Omega \subset\mathbb{R}^3$ initially spans a rectangular parameter domain $\tilde{\Omega}$. To represent complex topologies, such as a surface with holes, the concept of trimming is introduced.

The core idea is to define closed loops in the parameter space that partition the domain into an active part $\tilde{\Omega}_{a}$ and a void part $\tilde{\Omega}_{v}$. These loops consist of connected trimming curves $\tilde{\mathbf{C}}_{k,j}(\tilde{\xi})$, which are themselves represented as B-spline curves:
\begin{equation}
\tilde{\mathbf{C}}_{k,j}(\tilde{\xi}) = 
\left[\begin{array}{c}
\xi_{k,j}(\tilde{\xi}) \\
\eta_{k,j}(\tilde{\xi})
\end{array}\right] = 
\sum_{i=1}^{n} N_i(\tilde{\xi})\, \tilde{\mathbf{P}}_i, 
\quad \tilde{\mathbf{P}}_i \in \mathbb{R}^2,
\end{equation}
where $N_i$ is the $i$-th univariate basis function and $\tilde{\mathbf{P}}_i$ are the control points in parameter space. The orientation of these curves determines the active area, typically following the convention that the domain to the left of the loop remains active. The mapping of the active surface is denoted as $\mathbf{S}_{a}:\tilde{\Omega}_{a} \rightarrow \Omega_{a}$. Details regarding trimming and the element-wise integration of trimmed domains can be found in our previous work \cite{hollweck_LR_THB_2026} and the cited literature therein.
\\

\begin{figure}[h]
    \centering
    \includegraphics[width=\textwidth]{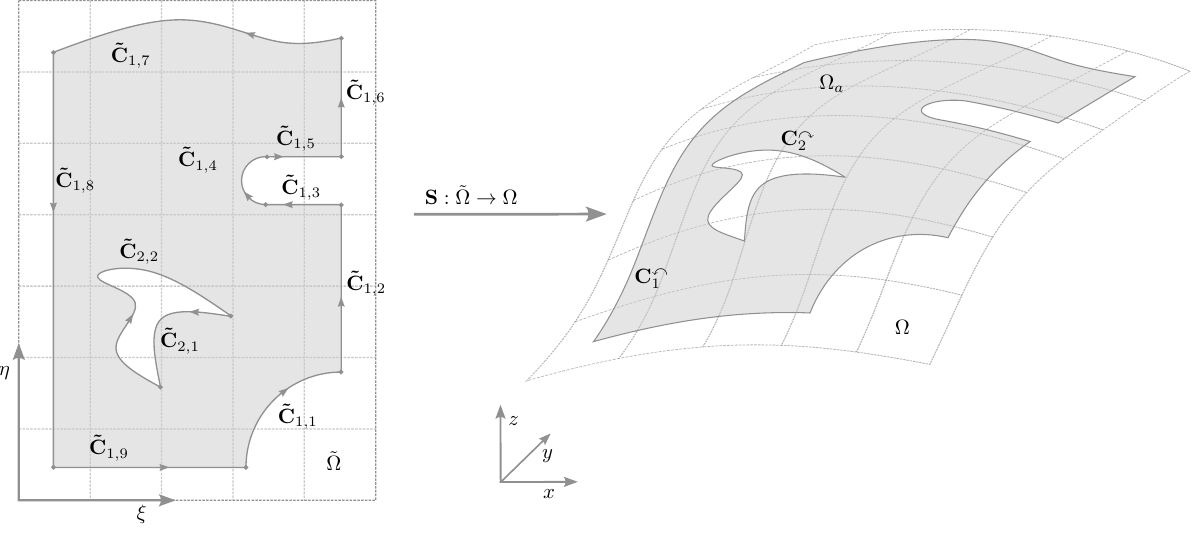}
    \caption{Mapping of a trimmed B-spline surface from parameter space to physical space. The geometry is defined by one outer loop $\mathbf{\tilde{C}}^{\curvearrowleft}_1$ and one inner loop $\mathbf{\tilde{C}}^{\curvearrowright}_2$. Figure adopted from \cite{hollweck_LR_THB_2026}.}
    \label{fig:sdsf}
\end{figure}

\subsection{Elastodynamic Model and Discretization}
\label{sec:elastodynamics}

\subsubsection*{Linear elastodynamics}

\rev{The linear elastodynamic model considered in the following provides a controlled setting for isolating discretization-induced spectral effects and, in particular, the upper spectral limit governing the explicit stability condition. It is not used to assess the accuracy of the corresponding strongly nonlinear dynamic response. The applicability of the resulting refinement strategy to such a nonlinear setting is assessed separately in the sheet metal forming benchmark of Sec.~\ref{sec:sec4}.}

We consider the linear elastodynamic initial boundary value problem on a domain 
$\Omega \subset \mathbb{R}^d$, $d \in \{1,2\}$, over a time interval $(0,T]$. 
Let $\mathbf{u} : \Omega \times (0,T] \rightarrow \mathbb{R}^d$ denote the displacement field. 
The strong form reads
\begin{equation}
\rho \ddot{\mathbf{u}} - \nabla \cdot \boldsymbol{\sigma} = \mathbf{f}
\quad \text{in } \Omega,
\label{eq:strong_dyn}
\end{equation}
where $\rho$ is the density, $\mathbf{f}$ the body force density, and 
$\boldsymbol{\sigma}$ the Cauchy stress tensor. 
The constitutive relation for linear elasticity is
\begin{equation}
\boldsymbol{\sigma} = \mathbf{C} : \boldsymbol{\varepsilon}(\mathbf{u}),
\qquad
\boldsymbol{\varepsilon}(\mathbf{u})
= \frac{1}{2}\left(\nabla \mathbf{u} + \nabla \mathbf{u}^{\mathrm{T}}\right).
\end{equation}
In one spatial dimension, Eq.~\eqref{eq:strong_dyn} reduces to
\begin{equation}
\rho \rev{A_{\mathrm{c}}} \ddot{u}
-
\frac{\partial}{\partial x}
\left(
E \rev{A_{\mathrm{c}}} \frac{\partial u}{\partial x}
\right)
=
f,
\label{eq:strong_1d}
\end{equation}
with cross-sectional area \rev{$A_{\mathrm{c}}$} and Young's modulus $E$.
The boundary is decomposed into disjoint parts 
$\Gamma_D$ and $\Gamma_N$ with prescribed displacement and traction,
\begin{equation}
\mathbf{u} = \mathbf{u}_D \text{ on } \Gamma_D,
\qquad
\rev{\boldsymbol{\sigma}\mathbf{n} = \mathbf{t}} \text{ on } \Gamma_N,
\end{equation}
and initial conditions
\begin{equation}
\mathbf{u}(\mathbf{x},0) = \mathbf{u}_0(\mathbf{x}),
\qquad
\dot{\mathbf{u}}(\mathbf{x},0) = \dot{\mathbf{u}}_0(\mathbf{x}).
\end{equation}

\subsubsection*{Weak Formulation}

Let
\[
\mathbf{V} = [H^1(\Omega)]^d
\]
denote the standard Sobolev space of square-integrable functions 
with square-integrable first derivatives.
Multiplying Eq.~\eqref{eq:strong_dyn} by a test function 
$\mathbf{v} \in \mathbf{V}$ and integrating by parts yields the weak form:

Find $\mathbf{u}:(0,T)\to\mathbf{V}$ such that for all $t\in(0,T]$
and for all $\mathbf{v} \in \mathbf{V}$,
\begin{equation}
m(\ddot{\mathbf{u}}(t),\mathbf{v})
+
a(\mathbf{u}(t),\mathbf{v})
=
l(\mathbf{v}),
\label{eq:weak_dyn}
\end{equation}
with
\begin{align}
m(\ddot{\mathbf{u}}(t),\mathbf{v})
&=
\int_{\Omega}
\rho \ddot{\mathbf{u}}(t) \cdot \mathbf{v}
\,\mathrm{d}\Omega,
\\
a(\mathbf{u}(t),\mathbf{v})
&=
\int_{\Omega}
\boldsymbol{\varepsilon}(\mathbf{u}(t))
:
\mathbf{C}
:
\boldsymbol{\varepsilon}(\mathbf{v})
\,\mathrm{d}\Omega,
\\
l(\mathbf{v})
&=
\int_{\Omega}
\mathbf{f} \cdot \mathbf{v}
\,\mathrm{d}\Omega
+
\int_{\Gamma_N}
\mathbf{t} \cdot \mathbf{v}
\,\mathrm{d}\Gamma.
\end{align}

\subsubsection*{Spline Discretization}

Let $\mathcal{M}$ be a locally refined Cartesian background mesh
on the parametric domain $\tilde{\Omega}$,
obtained either by LR- or THB-refinement.
Independently of the underlying specific formulation,
we denote the spline basis by
\[
\mathcal{T}
=
\{ T_A \}_{A\in\mathcal{I}},
\]
where $\mathcal{I}$ is the index set of active basis functions after trimming.
The symbol $T_A$ refers generically to a spline basis function,
either of LR- or THB-type.
In the trimmed setting, only basis functions whose support intersects
the active parametric domain $\tilde{\Omega}_a$ are retained, i.e.,
\[
\mathcal{I}
:=
\{\, A : \operatorname{supp}(T_A)\cap \tilde{\Omega}_a \neq \emptyset \,\}.
\]
With the geometric mapping $\mathbf{S}:\tilde{\Omega}\rightarrow\Omega$
and $\Omega_a=\mathbf{S}(\tilde{\Omega}_a)$, the discrete space reads
\[
\mathbf{V}_h(\Omega_a)
:=
\left\{
\mathbf{v}_h :
(\mathbf{v}_h \circ \mathbf{S})
=
\sum_{A\in\mathcal{I}} T_A \mathbf{v}_A,
\quad
\mathbf{v}_A \in \mathbb{R}^d
\right\}
\subset [H^1(\Omega_a)]^d.
\]
While LR-spline basis functions are tensor products on their support,
THB-splines are not globally tensor-product functions due to truncation.
Nevertheless, their restriction to any element $e\in\mathcal{M}$
remains a tensor-product polynomial,
\[
T_A|_e \in \mathcal{Q}^{p}(e),
\]
where $p$ denotes the spline degree.
The same conforming space is used for both trial and test functions.
Using the expansion
\[
\mathbf{u}_h(\mathbf{x},t)
=
\sum_{A\in\mathcal{I}} T_A(\mathbf{x})\,\mathbf{u}_A(t),
\]
where $\mathbf{u}_A(t)\in\mathbb{R}^d$ denote the displacement coefficients,
we collect all coefficients into the global vector
\[
\mathbf{u}(t)
=
\begin{bmatrix}
\mathbf{u}_1(t)\\
\mathbf{u}_2(t)\\
\vdots
\end{bmatrix}
\in\mathbb{R}^{d|\mathcal{I}|}.
\]

\begin{remark}
The above representation makes the separation of spatial and temporal
dependencies explicit: the spline basis functions $T_A$ describe the
spatial approximation, while the coefficients $\mathbf{u}_A(t)$ carry
the entire time dependence in a semidiscrete Galerkin sense.
Strictly speaking, $T_A$ are defined on the parametric domain
$\tilde{\Omega}$ and the physical representation reads
$T_A \circ \mathbf{S}^{-1}(\mathbf{x})$ with the geometric mapping
$\mathbf{S}:\tilde{\Omega}\to\Omega$.
\end{remark}

Let $A=1,\dots,|\mathcal{I}|$ denote the basis function index and 
$\alpha=1,\dots,d$ the spatial component.
We introduce the global degree-of-freedom index
\[
i = d(A-1)+\alpha,
\]
so that the global vector satisfies
\[
u_i = (\mathbf{u}_A)_\alpha.
\]
Insertion into the weak form yields the matrix system
\[
\mathbf{M}\ddot{\mathbf{u}}(t)
+
\mathbf{K}\mathbf{u}(t)
=
\mathbf{F}(t),
\]
where the matrix entries are given by
\[
M_{ij}
=
m\!\left(T_A \mathbf{e}_\alpha,\; T_B \mathbf{e}_\beta\right),
\qquad
K_{ij}
=
a\!\left(T_A \mathbf{e}_\alpha,\; T_B \mathbf{e}_\beta\right),
\]
with $i=d(A-1)+\alpha$ and $j=d(B-1)+\beta$.

\begin{remark}
In Isogeometric Analysis, Dirichlet boundary conditions cannot in general be imposed strongly, 
since spline basis functions are typically non-interpolatory.
Therefore, essential boundary conditions are commonly enforced weakly, 
for instance by means of penalty formulations, Nitsche's method, or Lagrange multipliers \cite{leidinger_explicit_2019, Coupling, Lagrangemult, nitsche}.
\rev{For the spectral investigations in this work, the undamped free-vibration problem is considered without essential boundary conditions or external loads. This free-free configuration isolates spectral effects introduced by the discretization, trimming, and local refinement from additional effects associated with boundary-condition enforcement.}
\end{remark}

\subsection{Numerical Model Setup}
\label{sec:models}

For the numerical investigations, we employ 1D bar and 2D membrane models.
Units are omitted for simplicity. All integrals are formulated in the
parametric domain and mapped to the physical domain via the geometric
Jacobian $\mathbf{J}=\partial \boldsymbol{x}/\partial \boldsymbol{\xi}$.

\subsubsection*{1D Bar Model}

The bar is discretized into $n_e$ elements using equidistant open knot vectors.
We assume a linear elastic material with Young's modulus $E=1$,
density $\rho=1$, and cross-sectional area \rev{$A_{\mathrm{c}}=1$}.
\rev{A diagonal mass matrix obtained by standard row-sum lumping is employed.}
The diagonal mass entries read
\begin{equation}
M_{ii}
=
\rho \rev{A_{\mathrm{c}}}
\int_{\tilde{\Omega}_a}
T_i
\, \left|\frac{\partial x}{\partial \xi}\right|
\, \mathrm{d}\xi,
\end{equation}
while the stiffness matrix entries are given by
\begin{equation}
K_{ij}
=
E \rev{A_{\mathrm{c}}}
\int_{\tilde{\Omega}_a}
\frac{\partial T_i}{\partial x}
\frac{\partial T_j}{\partial x}
\left|\frac{\partial x}{\partial \xi}\right|
\, \mathrm{d}\xi.
\end{equation}
For validation, numerical eigenfrequencies are compared with the analytical
solution of a free-free bar of length $L$ given by
\begin{equation}
\omega_n
=
\frac{n\pi}{L}
\sqrt{\frac{E}{\rho}},
\qquad
n=1,2,3,\dots.
\label{eq:free_rod}
\end{equation}

\subsubsection*{2D Membrane Model}

The membrane is defined on a square domain ($L_1=L_2$)
with thickness $t=1$ and Poisson ratio $\nu=0.3$.
\rev{Again, a diagonal mass matrix obtained by standard row-sum lumping is employed.}
Each control point carries identical mass in both displacement directions.
The lumped control point mass associated with basis function $T_A$ is
\begin{equation}
M_A
=
\rho t
\int_{\tilde{\Omega}_a}
T_A
\, \left|\det \mathbf{J}\right|
\, \mathrm{d}\tilde{\Omega}_a,
\end{equation}
where $\mathbf{J}=\partial \boldsymbol{x}/\partial \boldsymbol{\xi}$
denotes the geometric Jacobian.
Since two displacement components are present, the diagonal mass entries satisfy
\[
M_{ii} = M_A,
\qquad
i = 2(A-1)+\alpha,
\quad
\alpha=1,2.
\]
The element stiffness matrix is obtained from
\begin{equation}
\mathbf{K}^e
=
t
\int_{\tilde{\Omega}_a^e}
\mathbf{B}^{e\mathrm T}
\mathbf{C}
\mathbf{B}^e
\, \left|\det \mathbf{J}\right|
\, \mathrm{d}\tilde{\Omega}_a,
\end{equation}
where Voigt notation is employed. Accordingly,
$\mathbf{C}$ is given in its $3\times3$ plane-stress form.
The strain-displacement matrix reads
\[
\mathbf{B}^e
=
\begin{bmatrix}
\partial_x T_1^e & 0 & \cdots & \partial_x T_k^e & 0 \\
0 & \partial_y T_1^e & \cdots & 0 & \partial_y T_k^e \\
\partial_y T_1^e & \partial_x T_1^e & \cdots &
\partial_y T_k^e & \partial_x T_k^e
\end{bmatrix},
\]
where $k$ denotes the number of basis functions
with support on element $e$.
Plane stress is assumed with constitutive matrix
\[
\mathbf{C}
=
\frac{E}{1-\nu^2}
\begin{bmatrix}
1 & \nu & 0 \\
\nu & 1 & 0 \\
0 & 0 & \frac{1-\nu}{2}
\end{bmatrix}.
\]
In the numerical experiments presented below,
parametric and physical coordinates coincide,
so that $\left|\det \mathbf{J}\right|=1$.

\subsection{Stability of Explicit Time Integration}
\label{sec:stability}

Explicit time integration is conditionally stable, requiring a time step $\Delta t$ below a specific critical threshold $\Delta t_{\text{crit}}$. This approach is particularly advantageous for short-term dynamic simulations involving significant nonlinearities, such as crashworthiness or metal-forming analyses, where fine temporal resolution is physically required. Unlike implicit schemes, explicit methods avoid solving large systems of equations and circumvent convergence issues related to structural instabilities like buckling or wrinkling \cite{benson_large_2011}, provided the application of a lumped mass matrix.

The stability limit is dictated by the maximum system eigenfrequency $\omega_{\max}$, derived from the generalized eigenvalue problem:
\begin{equation}
    (\mathbf{K} - \omega_i^2 \mathbf{M}) \mathbf{q}_i = \mathbf{0},
\end{equation}
where $\mathbf{M}$ denotes the (usually lumped) mass matrix, $\mathbf{K}$ the stiffness matrix, while $\omega_i^2$ and $\mathbf{q}_i$ denote the $i$-th eigenvalue and eigenvector, respectively. For the central difference method and an undamped system, the stable time step must satisfy the Courant-Friedrichs-Lewy (CFL) condition:
\begin{equation}
    \Delta t \leq \Delta t_{\text{crit}} = \frac{2}{\omega_{\max}}.
    \label{eq:inequality_timestep}
\end{equation}

\rev{Since determining $\omega_{\max}$ for large systems at every time step is computationally expensive,} efficient estimation techniques are often employed, such as nodal or element-based estimates \cite{HUGHES_Implicit, BELYTSCHKO_Stability}, where the estimate should always be conservative but as sharp as possible.

\rev{A common conservative approach employs the classical bound relating the maximum system eigenfrequency to the largest element eigenfrequency \cite{IronsTreharne1971,Hughes2012}:}
\begin{equation}
\omega_{\max}
\leq
\max_{e} (\omega_{\max}^{e})
=
\omega_{\max}^{\text{E}}.
\end{equation}

Instead of solving the eigenvalue problem at the element level, the element geometry often serves as a heuristic for the maximum element eigenfrequency.
\rev{Such estimates typically rely on a characteristic element length $\ell_c^{e}$ and the material's elastic wave speed $c = \sqrt{\tilde{E}/\rho}$, where $\tilde{E}$ is the appropriate elastic modulus, for example for plane stress or plane strain. Corresponding element- and nodal-based time-step estimates for smooth B-spline and NURBS discretizations, including their dependence on polynomial degree, continuity, and boundary elements, were investigated in detail by Adam et al.~\cite{ADAM2015581}.}
For an undamped system, the maximum frequency is related to these physical quantities via \cite{Hartmann_stable}:
\begin{equation}
    \omega_{\max} \approx \frac{2 c}{\min_{e}({\ell_c^{e}})}.
\end{equation}

Substituting this into Eq.~\eqref{eq:inequality_timestep} yields the classic CFL condition $\Delta t \leq \min_{e}(\ell_c^{e}) / c$. While $\ell_c^{e}$ is a heuristic value depending on the element geometry, for example the minimum edge length or diagonal, it provides a computationally inexpensive alternative to spectral analysis. In industrial setups, a safety factor $f_s$ (typically $f_s = 0.9$) is employed to account for nonlinearities. Consequently, the stable time step is determined by:
\begin{equation}
    \Delta t \leq  f_s \frac{\min_{e}({\ell_c^{e}})}{c}.
\end{equation}

While this method originates from standard FEM, for trimmed isogeometric meshes it is not obvious what the characteristic element length $\ell_c^{e}$ of a trimmed element actually is. However, a substituted characteristic element length $\ell_c^{e}$ based on the Jacobian can be found in \cite{leidinger_Diss, Hartmann_stable} that also accounts for different polynomial degrees and higher continuity. An adjusted time step estimator thereof for IGA is currently used in LS-DYNA \cite{hallquist2006lsdyna}.

Lastly, the Gershgorin circle theorem \cite{gerschgorin} provides a conservative upper bound for the maximum system eigenfrequency. For a diagonal mass matrix of size $n \times n$, where $n = d|\mathcal{I}|$ denotes the total number of degrees of freedom, the eigenvalues $\omega_i^2$ are bounded by:
\begin{equation}
    \omega_{\max}^2
    \leq
    \max_i
    \left(
    \frac{1}{M_{ii}}
    \sum_{j=1}^{n}
    |K_{ij}|
    \right)
    =
    \bigl(\omega^{\mathrm{G}}_{\max}\bigr)^2.
\label{eq:gersh}
\end{equation}

Alternatively, a modified Gershgorin bound can be derived as reported in \cite{Bioli2025} from the $\infty$-matrix norm of the symmetrically scaled stiffness matrix:
\begin{equation}
    \omega_{\max}^2
    \leq
    \left\|
    \mathbf{M}^{-1/2}
    \mathbf{K}
    \mathbf{M}^{-1/2}
    \right\|_\infty
    =
    \max_i
    \left(
    \sum_{j=1}^{n}
    \frac{|K_{ij}|}{\sqrt{M_{ii}M_{jj}}}
    \right)
    =
    \bigl(\omega^{\mathrm{G,mod}}_{\max}\bigr)^2.
\label{eq:mod_gersh}
\end{equation}
It provides a more precise bound for the maximum eigenfrequency in the presence of non-uniform mass distributions, such as those caused by trimming.

\section{Numerical Investigation of Trimming-Induced High-Frequency Outliers}
\label{sec:sec3}
This section investigates the spectral properties of trimmed B-spline patches of maximum continuity $C^{p-1}$ under various local refinement strategies. The primary objective is to analyze how trimming-induced support reduction and the spatial distribution of refined basis functions affect the maximum patch eigenfrequency $\omega_{\text{max}}$ and, consequently, the critical time step of explicit time integration.

\rev{The investigations proceed from controlled one-dimensional configurations to two-dimensional trimmed geometries. Particular attention is paid to basis functions with reduced active support adjacent to the trimming boundary and to their interaction with local refinement. Knot-exact and arbitrary trimming are considered as different configurations of the same underlying support-reduction mechanism. Based on these observations, the Boundary-Level-Constrained Refinement (BLCR) strategy is introduced to prevent refined trimmed basis functions from governing the upper spectral limit.}

\rev{The primary quantity of interest is the maximum eigenfrequency governing the critical time step. Nevertheless, changes in the lower part of the spectrum caused by the refinement constraint are also examined where relevant, since these changes may affect approximation accuracy. The present analysis does not attempt to resolve the general spectral-accuracy limitations associated with row-sum mass lumping. All spectral investigations are performed for undamped free-vibration problems without essential boundary conditions.}

\subsection{Investigations in 1D}

We begin our numerical study with one-dimensional investigations to isolate the fundamental effects of trimming and local refinement on the maximum eigenfrequency. The 1D setting serves as a benchmark, allowing for a highly controlled environment where boundary effects can be studied analytically and numerically. 

\subsubsection{Refinement and Trimming Setups}
\label{sec:refinement_strategies}

\begin{figure}[H]
    \centering
    \renewcommand{\arraystretch}{1.4}

    \newcommand{\LeftBox}[1]{\makebox[0.375\textwidth][c]{#1}}
    \newcommand{\RightBox}[1]{\makebox[0.625\textwidth][c]{#1}}

    \begin{tabular}{@{} c c @{}}
        \LeftBox{
            \subfloat[][B-spline, untrimmed, $\omega_{\text{max}}=1.16$]{\includegraphics[scale=0.5]{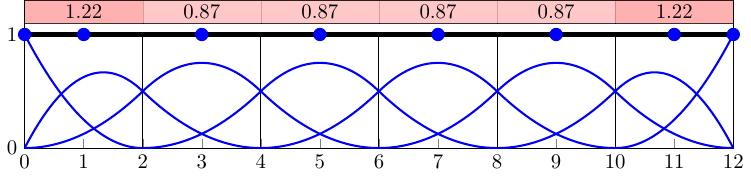}}
        } &
        \RightBox{
            \subfloat[][B-spline, boundary-trimmed, $\omega_{\text{max}}=0.75$]{\includegraphics[scale=0.5]{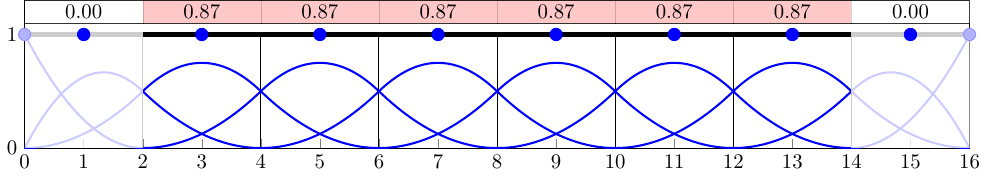}}
        } \\[4ex]
        \LeftBox{
            \subfloat[][globally refined B-spline, untrimmed, $\omega_{\text{max}}=2.32$]{\includegraphics[scale=0.5]{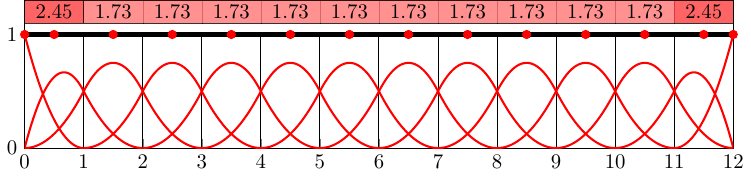}}
        } &
        \RightBox{
            \subfloat[][globally refined B-spline, boundary-trimmed, $\omega_{\text{max}}=1.49$]{\includegraphics[scale=0.5]{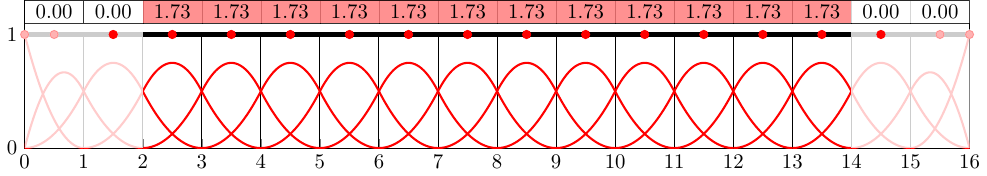}}
        } \\[4ex]
        \LeftBox{
            \subfloat[][BLCR with LR-spline, untrimmed, $\omega_{\text{max}}=1.23$]{\includegraphics[scale=0.5]{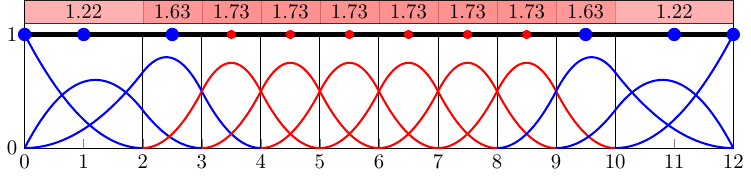}}
        } &
        \RightBox{
            \subfloat[][BLCR with LR-spline, boundary-trimmed, $\omega_{\text{max}}=1.27$]{\includegraphics[scale=0.5]{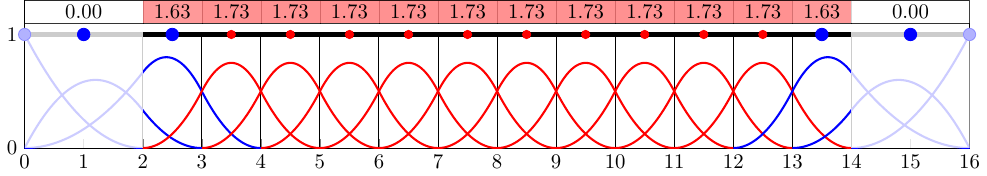}}
        } \\[4ex]
        \LeftBox{
            \subfloat[][BLCR with THB-spline, untrimmed, $\omega_{\text{max}}=1.23$]{\includegraphics[scale=0.5]{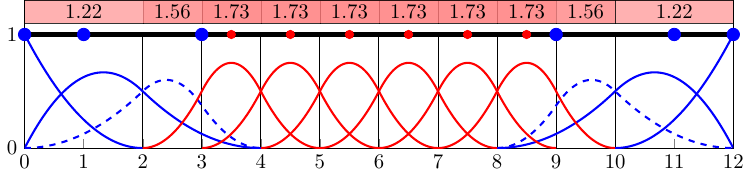}}
        } &
        \RightBox{
            \subfloat[][BLCR with THB-spline, boundary-trimmed, $\omega_{\text{max}}=1.23$]{\includegraphics[scale=0.5]{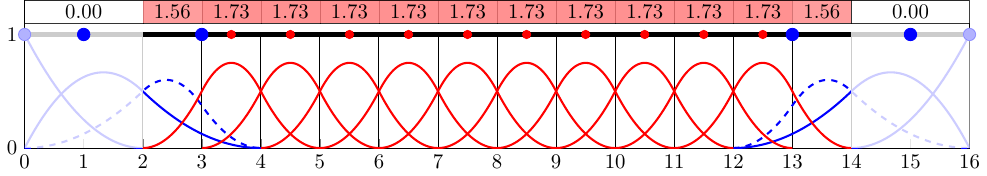}}
        } \\
    \end{tabular}

     \caption{One-dimensional investigations for $p=2$: The maximum element eigenfrequencies $\omega_{\text{max}}^e$ are indicated above each element. Inactive control points are masked for clarity. Dotted lines represent truncated basis functions. Blue denotes coarse basis functions, while red represents refined basis functions.}
    \label{fig:refinement_p2}
\end{figure}

\begin{figure}[H]
    \centering
    \renewcommand{\arraystretch}{1.4}

    \newcommand{\LeftBox}[1]{\makebox[0.375\textwidth][c]{#1}}
    \newcommand{\RightBox}[1]{\makebox[0.625\textwidth][c]{#1}}

    \begin{tabular}{@{} c @{\hspace{-0.09cm}} c @{}}
        \LeftBox{
            \subfloat[][B-spline, untrimmed, $\omega_{\text{max}}=1.53$]{\includegraphics[scale=0.49]{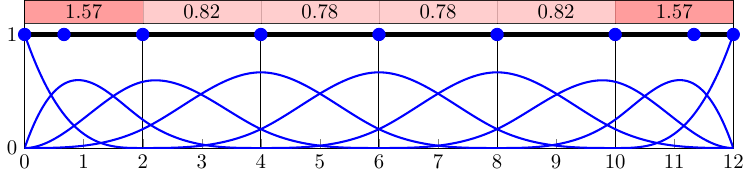}}
        } &
        \RightBox{
            \subfloat[][B-spline, boundary-trimmed, $\omega_{\text{max}}=0.66$]{\includegraphics[scale=0.49]{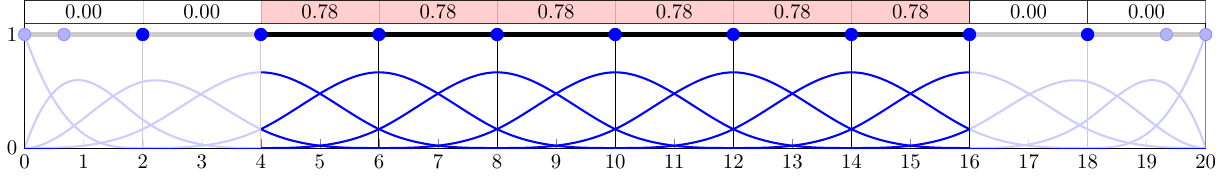}}
        } \\[4ex]
        \LeftBox{
            \subfloat[][globally refined B-spline, untrimmed, $\omega_{\text{max}}=3.06$]{\includegraphics[scale=0.49]{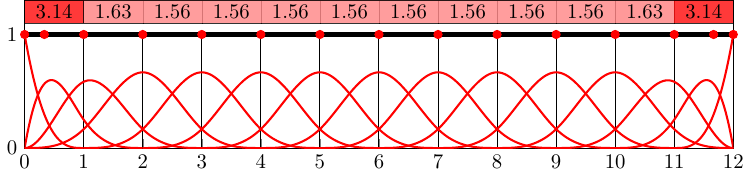}}
        } &
        \RightBox{
            \subfloat[][globally refined B-spline, boundary-trimmed, $\omega_{\text{max}}=1.32$]{\includegraphics[scale=0.49]{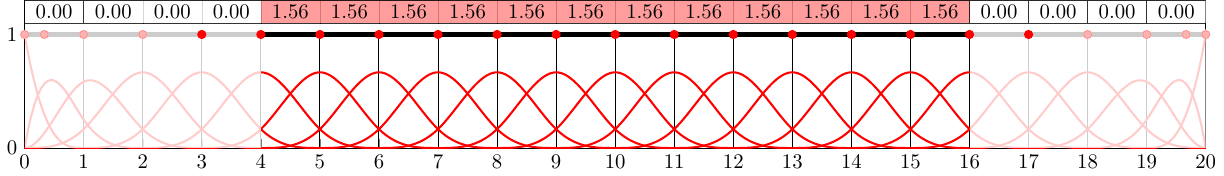}}
        } \\[4ex]
        \LeftBox{
            \subfloat[][BLCR with LR-spline, untrimmed, $\omega_{\text{max}}=1.57$]{\includegraphics[scale=0.49]{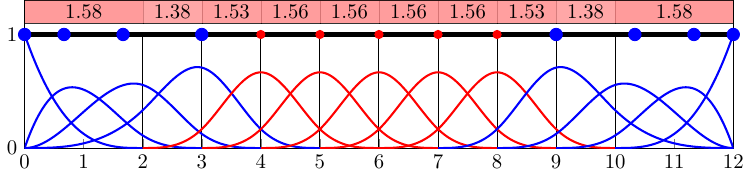}}
        } &
        \RightBox{
            \subfloat[][BLCR with LR-spline, boundary-trimmed, $\omega_{\text{max}}=1.08$]{\includegraphics[scale=0.49]{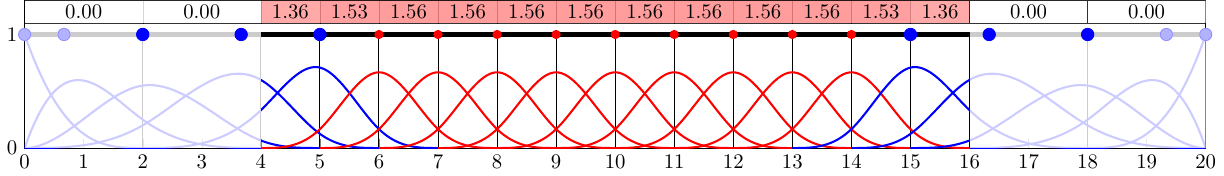}}
        } \\[4ex]
        \LeftBox{
            \subfloat[][BLCR with THB-spline, untrimmed, $\omega_{\text{max}}=1.53$]{\includegraphics[scale=0.49]{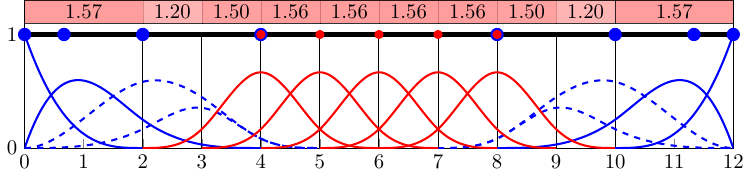}}
        } &
        \RightBox{
            \subfloat[][BLCR with THB-spline, boundary-trimmed, $\omega_{\text{max}}=1.02$]{\includegraphics[scale=0.49]{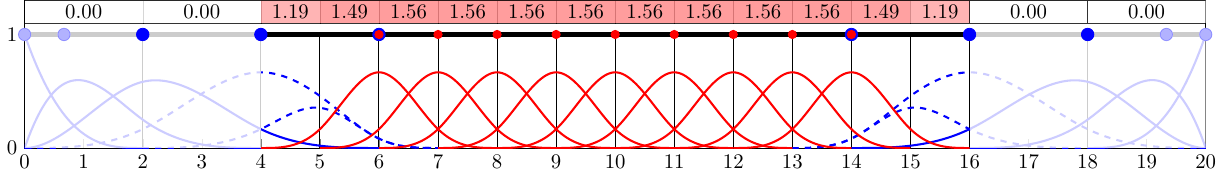}}
        } \\
    \end{tabular}

     \caption{One-dimensional investigations for $p=3$: The maximum element eigenfrequencies $\omega_{\text{max}}^e$ are indicated above each element. Inactive control points are masked for clarity. Dotted lines represent truncated basis functions. Blue denotes coarse basis functions, while red represents refined basis functions.}
    \label{fig:refinement_p3}
\end{figure}

We begin our investigations using a 1D rod model (\(L=12\)), where we study how global and local refinement strategies, combined with the established boundary-trimming technique, affect \( \omega_{\text{max}} \) for degrees \(p=2\) (Fig.~\ref{fig:refinement_p2}) and \(p=3\) (Fig.~\ref{fig:refinement_p3}). The numbers above elements indicate the maximum element eigenfrequency \( \omega_{\text{max}}^e \).

\paragraph{Standard B-spline (Fig.~\ref{fig:refinement_p2}a,b)}
The uniform B-spline discretization with six elements in Fig.~\ref{fig:refinement_p2}a yields a patch eigenfrequency of \( \omega_{\text{max}} = 1.16 \). A pronounced boundary effect is visible, as the first and last elements exhibit significantly higher \( \omega_{\text{max}}^e \) due to the altered shape of the basis functions at the patch edges. Applying the standard remedy by extending the parametric domain and trimming off the boundary \(p-1\) knot spans \cite{leidinger_Diss, hollweck_LR_THB_2026, Reali_explicit} results in a configuration with only interior, uniform elements, reducing the global \( \omega_{\text{max}} \) to 0.75 (Fig.~\ref{fig:refinement_p2}b).

\paragraph{Global Refinement (Fig.~\ref{fig:refinement_p2}c,d)}
Global \(h\)-refinement (knot insertion) of both untrimmed (Fig.~\ref{fig:refinement_p2}c) and trimmed (Fig.~\ref{fig:refinement_p2}d) configurations halves the characteristic element length and, as expected, scales the largest eigenfrequency by a factor of two. \rev{In the knot-exact trimmed case, all active elements share an identical \( \omega_{\text{max}}^e \). This element-wise indicator therefore no longer reveals the contribution of reduced-support basis functions adjacent to the trimming boundary, which becomes visible only in the assembled patch spectrum considered below.}

\paragraph{Interior Local Refinement (Fig.~\ref{fig:refinement_p2}e--h)}
\label{par:local_refinement}
Crucially, local refinement enables a further reduction of \( \omega_{\text{max}} \) compared to Fig.~\ref{fig:refinement_p2}c,d. Starting from the coarse meshes in Fig.~\ref{fig:refinement_p2}a (untrimmed) and Fig.~\ref{fig:refinement_p2}b (trimmed), we refine all basis functions except those directly associated with the physical boundary (for the trimmed case in Fig.~\ref{fig:refinement_p2}f,h) or the first and last functions (for the untrimmed case in Fig.~\ref{fig:refinement_p2}e,g), using both LR- and THB-splines.

\begin{itemize}
    \item For the \textbf{untrimmed} configurations (Fig.~\ref{fig:refinement_p2}e,g), refining the interior elevates the element eigenfrequencies in that region (e.g.\ \( \omega_{\text{max}}^e = 1.73 \)). These now surpass those of the intentionally unrefined boundary elements (\( \omega_{\text{max}}^e = 1.22 \)). This strategy reduces the patch \( \omega_{\text{max}} \) to 1.23, which is lower than both the globally refined and the standard trimmed configurations (Fig.~\ref{fig:refinement_p2}c,d).
    
    \item For the \textbf{trimmed} configurations (Fig.~\ref{fig:refinement_p2}f,h), a subtle difference appears: the LR-spline version (\( \omega_{\text{max}} = 1.27 \)) exhibits a slightly higher maximum than the THB-spline variant (\( \omega_{\text{max}} = 1.23 \)). As shown in Fig.~\ref{fig:refinement_p2}h, the coarse (solid and dotted blue) THB basis functions maintain support over two active elements, whereas the corresponding LR-spline functions in Fig.~\ref{fig:refinement_p2}f have a more restricted support. In Sec.~\ref{subsub:Rayleigh}, we demonstrate on a nodal basis that increasing the support of a trimmed basis function effectively reduces its eigenfrequency, which explains the observed spectral advantage of THB-splines in this configuration.
\end{itemize}

\paragraph{Behavior for $p=3$ (Fig.~\ref{fig:refinement_p3})}
The qualitative observations for \(p=3\) corroborate the trends. A key difference is that for the untrimmed, locally refined case (Fig.~\ref{fig:refinement_p3}e,g), the maximum eigenfrequency remains substantially higher (\(\omega_{\text{max}} = 1.57\) for LR-splines and \(\omega_{\text{max}} = 1.53\) for THB-splines) than in the trimmed counterpart (Fig.~\ref{fig:refinement_p3}f,h) (\(\omega_{\text{max}} = 1.08\) for LR-splines and \(\omega_{\text{max}} = 1.02\) for THB-splines). This suggests that for higher degrees, a single level of interior refinement may be insufficient to fully mitigate the original boundary effect, potentially requiring a multi-level strategy for untrimmed setups.

\paragraph{Key Finding and Interpretation}

\rev{The preceding investigations reproduce and extend several observations reported by Adam et al.~\cite{ADAM2015581}. Adam et al. identified the restrictive boundary effect of clamped (open-knot) B-splines and investigated two alternative remedies for this same limitation. The first is to increase the size of the boundary elements until the critical time step is governed by interior functions. The second is to unclamp the knot vector, thereby removing the particularly restrictive exterior boundary functions associated with the clamped representation. For the knot-exact configurations considered here, trimming away the exterior open-knot spans is operationally equivalent to this unclamping strategy and yields the corresponding active discrete representation.}

\rev{The present results show how these two principles can be combined. Removing the exterior open-knot functions first suppresses the dominant clamped-boundary penalty. However, basis functions adjacent to the resulting trimming boundary still possess reduced active support and can remain relevant for the maximum eigenfrequency. If these functions are refined together with the interior, they can become the limiting spectral contributors. BLCR therefore combines the unclamping-equivalent removal of the exterior boundary functions with the boundary-coarsening principle investigated by Adam et al., applying the latter to the basis functions that remain after the exterior open-knot functions have already been removed.}

\rev{This distinction is also visible in the untrimmed boundary-fitted configurations. For $p=2$, a difference of one refinement level is sufficient in the investigated cases to remove the boundary-induced limitation. For $p=3$, the same one-level difference in Figs.~\ref{fig:refinement_p3}e,g is no longer sufficient while the particularly restrictive clamped boundary functions remain active. In contrast, after these exterior functions are removed by knot-exact trimming, the additional relative coarsening in Figs.~\ref{fig:refinement_p3}f,h is sufficient in the investigated configurations. This stronger degree dependence of the clamped boundary contribution is consistent with the observations of Adam et al.}

\rev{The additional observation of the present work is that, even after the particularly restrictive exterior open-knot functions have been removed, reduced-support basis functions adjacent to a knot-exact trimming boundary can still govern the maximum eigenfrequency, as indicated by Figs.~\ref{fig:refinement_p2}d,f,h and \ref{fig:refinement_p3}d,f,h. These remaining basis functions were not isolated as the governing global spectral contributors in Adam et al., nor was unclamping combined with relative boundary coarsening.}

\rev{This remaining contribution is not apparent from the maximum element eigenfrequencies alone. In the globally refined boundary-trimmed configurations of Figs.~\ref{fig:refinement_p2}d and \ref{fig:refinement_p3}d, the active elements exhibit identical maximum element eigenfrequencies, yet the assembled system still contains an elevated maximum eigenfrequency associated with basis functions of reduced active support. Thus, the present study does not introduce a separate physical boundary mechanism. Instead, it identifies the remaining spectral driver after removal of the classical open-knot boundary functions and shows how removal of these functions can be combined with relative boundary coarsening to prevent trimmed basis functions from governing the upper spectral limit in the investigated configurations.}

\subsubsection{Boundary-Level-Constrained Refinement (BLCR)}
\label{sub:Def_BLCR}

The preceding observations motivate the Boundary-Level-Constrained Refinement (BLCR) strategy. 
Let $\tilde{\Omega}$ denote the parametric domain and $\tilde{\Omega}_a \subset \tilde{\Omega}$ the active parametric domain as introduced in Sec.~\ref{sec:trimmedsurfaces}. Let 
\begin{equation}
    \mathcal{T} = \{T_A\}_{A \in \mathcal{I}}
\end{equation}
denote the set of active spline basis functions after trimming. We assume that the outer $p-1$ rows and columns of elements in the initial tensor-product mesh are fully trimmed. This removes the boundary effects inherent to the original domain boundaries, specifically those regions where the basis functions exhibit reduced continuity due to the open knot vector. Consequently, all remaining active basis functions in $\mathcal{T}$ possess maximal continuity $C^{p-1}$ throughout the active domain $\tilde{\Omega}_a$. This condition is typically ensured during a preprocessing step. Since the basis functions are defined on the parametric domain,
their supports satisfy
\[
\operatorname{supp}(T_A) \subset \tilde{\Omega}.
\]
We define the sets of interior and cut basis functions by
\[
\mathcal{T}_{\mathrm{int}}
:=
\left\{
T_A \in \mathcal{T}
\;:\;
\operatorname{supp}(T_A) \subset \tilde{\Omega}_a
\right\},
\qquad
\mathcal{T}_{\mathrm{cut}}
:=
\mathcal{T} \setminus \mathcal{T}_{\mathrm{int}}.
\]

Let $\ell(T_A)$ denote the refinement level of basis function $T_A$.
The BLCR strategy controls refinement such that the smallest
parametric length scales occur exclusively in the interior
of $\tilde{\Omega}_a$. More precisely, refinement is performed such that the level constraint
\[
\max_{T_A \in \mathcal{T}_{\mathrm{cut}}} \ell(T_A)
\;<\;
\max_{T_B \in \mathcal{T}_{\mathrm{int}}} \ell(T_B)
\]
is satisfied at every refinement step.
Hence, basis functions whose support intersects the trimmed boundary
may in principle be refined, but their refinement level must remain
strictly below the maximum level attained in the interior. This explains the name \emph{Boundary-Level-Constrained Refinement}: the maximum refinement level of basis functions associated with the trimming boundary is constrained relative to that of interior basis functions.

In the special case of a single refinement step, this condition implies
that only interior basis functions are refined, while trimmed basis
functions remain at their initial level.
For multiple refinement steps, trimmed basis functions may be refined
provided their level remains strictly lower than the maximal interior level.

\rev{For all trimmed configurations investigated in this work, maintaining the trimmed basis functions one refinement level below the maximally refined interior is sufficient to ensure that the upper spectral limit is governed by untrimmed interior functions. The corresponding separation of the frequency bounds is analyzed in Sec.~\ref{subsub:Rayleigh}. This observation should not be interpreted as a universal guarantee for arbitrary polynomial degrees, trimming configurations, or refinement hierarchies.}

\begin{remark}
Since trimming enables the representation of arbitrarily complex geometries, the BLCR strategy is formulated for trimmed patches with void boundary elements. Boundary-fitted discretizations with open knot vectors were nevertheless also considered in the preceding study. \rev{In this setting, the refinement strategy must instead keep the active boundary functions sufficiently coarse. As observed for $p=3$, a one-level difference is not sufficient for the boundary-fitted configuration considered above. In contrast, for the trimmed configurations investigated in this work, a one-level difference is sufficient.}
\end{remark}

\subsubsection{Eigenvalue Spectrum under Refinement and Trimming}

\begin{figure}[H]
    \centering
    \renewcommand{\arraystretch}{1.4}

    \newcommand{\LeftBox}[1]{\makebox[0.481\textwidth][c]{#1}}
    \newcommand{\RightBox}[1]{\makebox[0.519\textwidth][c]{#1}}

    \begin{tabular}{@{} c c @{}}
        \LeftBox{
            \subfloat[][untrimmed]{\includegraphics[scale=0.65]{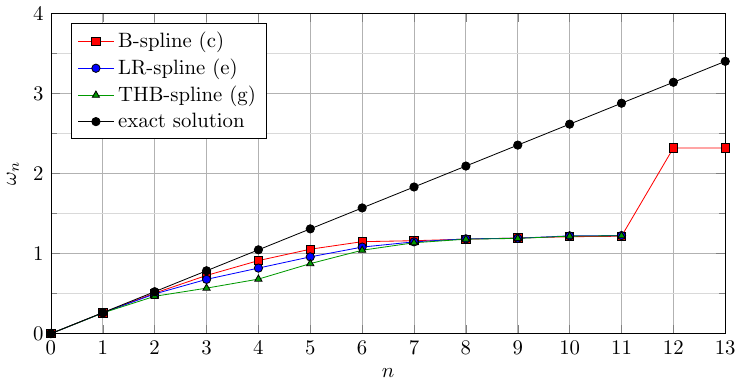}}
        } &
        \RightBox{
            \subfloat[][boundary-trimmed]{\includegraphics[scale=0.65]{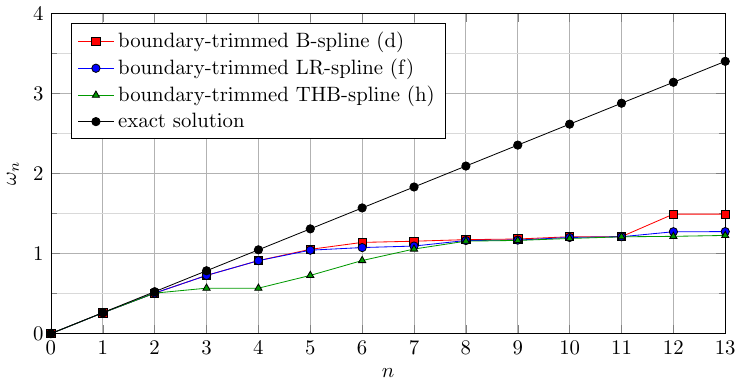}}
        }
    \end{tabular}

    \caption{Spectrum for $p=2$. Labels correspond to those in Fig.~\ref{fig:refinement_p2}.}
    \label{fig:spectrum_p2}
\end{figure}

\begin{figure}[H]
    \centering
    \renewcommand{\arraystretch}{1.4}

    \newcommand{\LeftBox}[1]{\makebox[0.481\textwidth][c]{#1}}
    \newcommand{\RightBox}[1]{\makebox[0.519\textwidth][c]{#1}}

    \begin{tabular}{@{} c c @{}}
        \LeftBox{
            \subfloat[][untrimmed]{\includegraphics[scale=0.65]{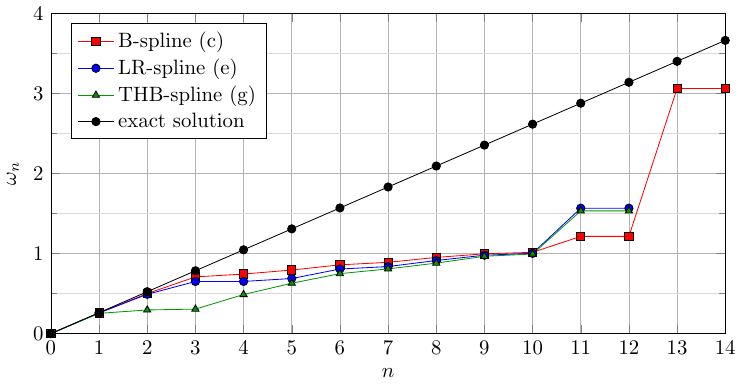}}
        } &
        \RightBox{
            \subfloat[][boundary-trimmed]{\includegraphics[scale=0.65]{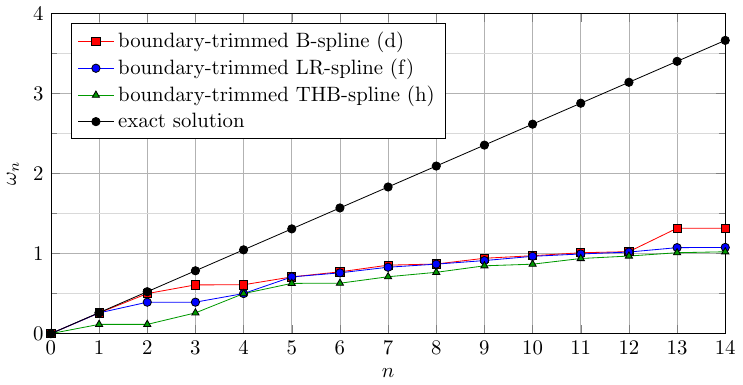}}
        }
    \end{tabular}

    \caption{Spectrum for $p=3$. Labels correspond to those in Fig.~\ref{fig:refinement_p3}.}
    \label{fig:spectrum_p3}
\end{figure}

\rev{The observation that local interior refinement reduces the maximum patch eigenfrequency relative to a globally refined trimmed patch indicates that reduced-support functions adjacent to the trimming boundary remain relevant even after the exterior open-knot boundary functions have been removed. Figs.~\ref{fig:spectrum_p2} and \ref{fig:spectrum_p3} present the spectra for the untrimmed (left) and trimmed (right) configurations for $p=2$ and $p=3$, respectively, for the previously investigated B-, LR-, and THB-spline discretizations. For reference, the analytical eigenfrequencies of the free-free 1D rod from Eq.~\ref{eq:free_rod} are also included.}

\begin{remark}
In the following plots, discrete values are indicated, e.g., by circles or squares. 
The markers are connected by thin solid lines solely to improve readability; 
these lines have no physical meaning.
\end{remark}

\rev{Row-sum mass lumping modifies the discrete spectrum relative to the consistent-mass formulation, and trimmed or immersed discretizations may exhibit spurious modes in the lower part of the spectrum \cite{Bioli2025, Voet_lumping_stabilization, Guarino2025}. We therefore consider Figs.~\ref{fig:spectrum_p2} and \ref{fig:spectrum_p3} not only with respect to $\omega_{\max}$, but also with respect to changes in the lower spectrum.}

\rev{For $p=2$, the standard untrimmed B-spline discretization in Fig.~\ref{fig:spectrum_p2}a exhibits the characteristic upper-spectrum contribution of the exterior open-knot basis functions. Trimming off the exterior boundary spans substantially reduces $\omega_{\max}$, but the last two frequencies remain elevated (Fig.~\ref{fig:spectrum_p2}b). Starting one refinement level coarser and refining the interior with LR- or THB-splines further reduces this contribution, which is the effect targeted by BLCR. Note that the standard B-spline discretization in Fig.~\ref{fig:spectrum_p2}a contains 14 DOFs, compared with 12 DOFs for the LR- and THB-spline variants. The indices shown in the figures start at zero.}

\rev{For $p=3$, boundary trimming again substantially reduces the upper spectral limit, while the locally refined boundary-fitted configurations retain a larger $\omega_{\max}$ than the knot-exact trimmed configurations in Fig.~\ref{fig:spectrum_p3}b. At the same time, the refinement constraint visibly affects individual frequencies in the lower part of the spectrum, most notably the first non-zero frequencies of the locally refined THB-spline discretizations in Fig.~\ref{fig:spectrum_p3}b. Remarkably, these deviations occur even though the trimming boundary coincides exactly with knot lines.}

\rev{Spurious low-frequency modes due to row-sum lumping have previously been reported for arbitrarily trimmed discretizations, particularly when the active support of individual basis functions becomes very small \cite{Bioli2025, Voet_lumping_stabilization, Guarino2025}. The present results show that noticeable deviations in the lower spectrum may also occur for knot-exact trimming, without partially cut elements or vanishingly small cut fractions. This suggests that reduced active support and the resulting changes in local stiffness-to-mass ratios are more general features relevant to these spectral changes than the presence of a small cut cell itself. A detailed analysis of the associated low-frequency modes is beyond the scope of the present work.}

\rev{These observations reveal an important trade-off of BLCR: reducing the upper spectral limit governing the critical time step can alter individual lower frequencies. BLCR should therefore be applied with particular care to low-frequency-dominated problems and to applications requiring accurate dynamics near the trimmed boundary, including knot-exact configurations. Its intended use is instead in settings where the critical explicit time step is the dominant efficiency constraint and the trimming boundary is not a primary region of approximation interest.}

\subsubsection{Refinement from left to right}
\begin{figure}[H]
    \centering
    \includegraphics[width=\textwidth]{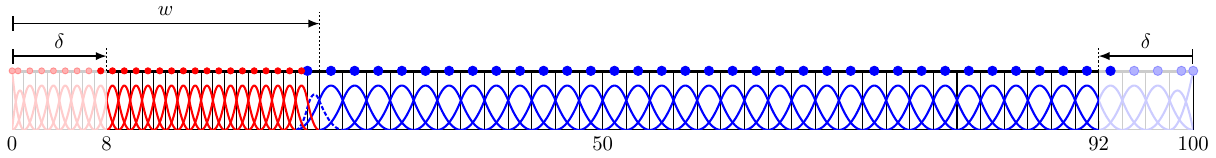}
    \caption{THB-spline, $p=2$, Refinement from left to right end. Trimming distance from both ends is fixed to $\delta=8$.}
    \label{fig:left2right_THB_p2}
\end{figure}

\begin{figure}[H]
    \centering
    \includegraphics[width=\textwidth]{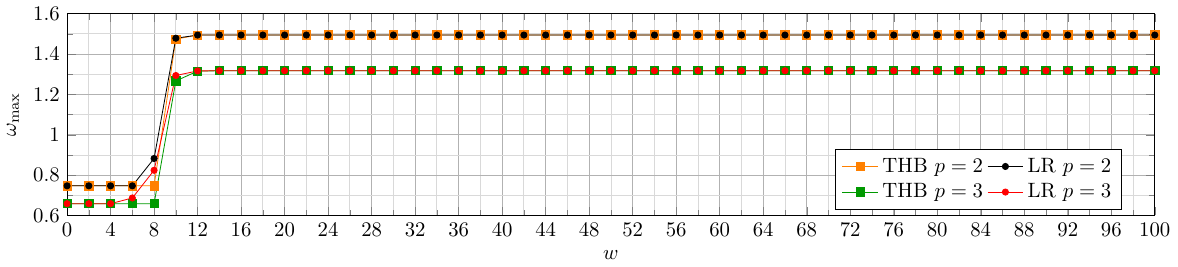}
    \caption{$\omega_{\text{max}}$ for THB- and LR-splines for $p=2$ and $p=3$.}
    \label{fig:omega_max_plot_left2right}
\end{figure}

We next isolate the effect of refining trimmed basis functions. Starting from a uniform $C^{p-1}$ B-spline discretization with 50 elements of width $2$, trimming is applied symmetrically from both ends with a fixed trimming distance $\delta=8$. Refinement is then activated from the left boundary over a growing discrete refinement width $w$. Fig.~\ref{fig:left2right_THB_p2} illustrates the setup for THB-splines, and Fig.~\ref{fig:omega_max_plot_left2right} reports the resulting maximum patch eigenfrequency $\omega_{\max}(w)$.

For $w<\delta$, refinement affects only basis functions whose support lies entirely in the void region. Consequently, the active discretization and its spectrum remain unchanged, and $\omega_{\max}$ is constant. Once $w\ge\delta$, refinement reaches basis functions with non-zero active support, i.e., trimmed basis functions adjacent to the trimming interface. At this point, $\omega_{\max}$ exhibits a sharp increase to the characteristic outlier level observed previously for boundary-trimmed refined configurations ($\omega_{\max}=1.49$ for $p=2$; $\omega_{\max}=1.32$ for $p=3$). Differences between LR- and THB-splines are small; near the activation threshold $w\approx\delta$, LR-splines typically yield slightly larger values of $\omega_{\max}$, consistent with the support differences discussed in Figs.~\ref{fig:refinement_p2}f,h and \ref{fig:refinement_p3}f,h.

\subsubsection*{Observations:}

\begin{itemize}

\item (i) Refining trimmed basis functions triggers a jump of $\omega_{\max}$ due to the appearance of high-frequency outliers at the upper end of the spectrum. 

\item (ii) After the first trimmed basis function is refined, further refinement of additional trimmed basis functions has only a weak influence on $\omega_{\max}$. In practice, refining even a single trimmed basis function already incurs most of the time-step penalty.

\end{itemize}

\subsubsection{Symmetric Refinement from Center to Boundary}

\begin{figure}[H]
    \centering
    \includegraphics[width=\textwidth]{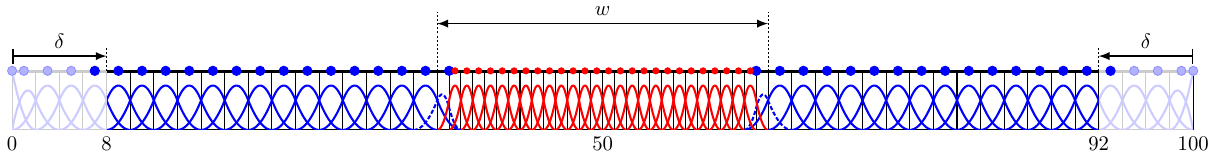}
    \caption{THB-spline, $p=2$, Symmetric refinement from center to end. Trimming distance from both ends is fixed to $\delta=8$.}
    \label{fig:center2end_THB_p2}
\end{figure}

\begin{figure}[H]
    \centering
    \includegraphics[width=\textwidth]{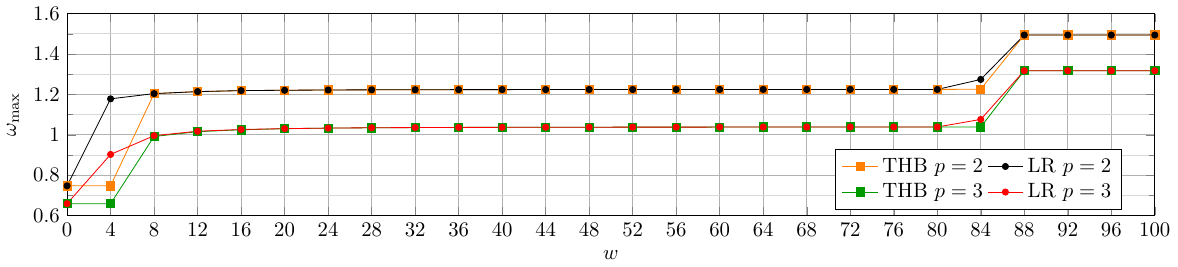}
    \caption{$\omega_{\text{max}}$ for THB- and LR-splines for $p=2$ and $p=3$.}
    \label{fig:omega_max_plot_center2end}
\end{figure}

We now consider the same setup as in the previous study, but now the refinement is initiated at the center of the patch and the refinement width is increased symmetrically towards the boundary, see Fig.~\ref{fig:center2end_THB_p2}. The differences between LR- and THB-splines for $w < 8$ in Fig.~\ref{fig:omega_max_plot_center2end} arise from the fact that, in the 1D case, the LR-spline refinement procedure already triggers a refinement, whereas for THB-splines the refinement width is still too small to fully contain at least one basis function. For $w \ge 8$, the same local refinement is applied to both LR- and THB-splines, and the resulting maximum patch eigenfrequencies are nearly identical.

When $w = 88$, trimmed basis functions are refined and the characteristic outlier frequencies reappear as in Fig.~\ref{fig:omega_max_plot_left2right} ($\omega_{\max}=1.49$ for $p=2$; $\omega_{\max}=1.32$ for $p=3$). For $w = 84$, the maximum patch eigenfrequency for LR-splines is again slightly higher than for THB-splines, since the knot insertion procedure can modify the shape of basis functions in coarse elements.

Overall, as long as refinement is restricted to the interior, the maximum patch eigenfrequency remains significantly lower than when trimmed boundary basis functions are refined. \rev{For the configurations investigated here, this shows that keeping the refined region separated from the trimming boundary prevents refined reduced-support functions from becoming the dominant contributors to $\omega_{\max}$. This observation forms the basis of BLCR.}

\subsubsection{Continuous Trimming from Boundary to Center}

\begin{figure}[H]
    \centering
    \includegraphics[width=0.4\textwidth]{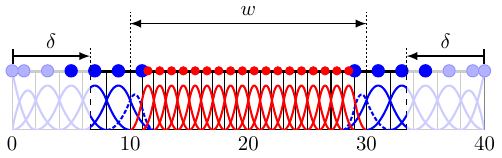}
    \caption{THB-spline, $p=2$, Fixed center refinement $w$. Trimming distance from both ends is denoted by $\delta$.}
    \label{fig:conti_trim_THB_p2}
\end{figure}

We conclude the 1D investigations by considering a setup in which local refinement is applied at the center with a fixed refinement width of $w = 20$, see Fig.~\ref{fig:conti_trim_THB_p2}. The boundary is trimmed by a distance $\delta$ towards the center, using both LR- and THB-splines.

For $p=2$, $\omega_{\max}$ remains constant until the trimming boundary reaches the refined region, after which it increases. For $p=3$, $\omega_{\max}$ initially decreases for very small trimming distances, becomes constant for $\delta>1$, and increases again once refined basis functions are affected by trimming. The increase starts slightly earlier for LR-splines than for THB-splines. When all coarse basis functions have been removed, the behavior approaches that of standard trimmed B-splines reported in \cite{leidinger_Diss, hollweck_LR_THB_2026, Reali_explicit}.

\begin{figure}[H]
    \centering
    \includegraphics[width=\textwidth]{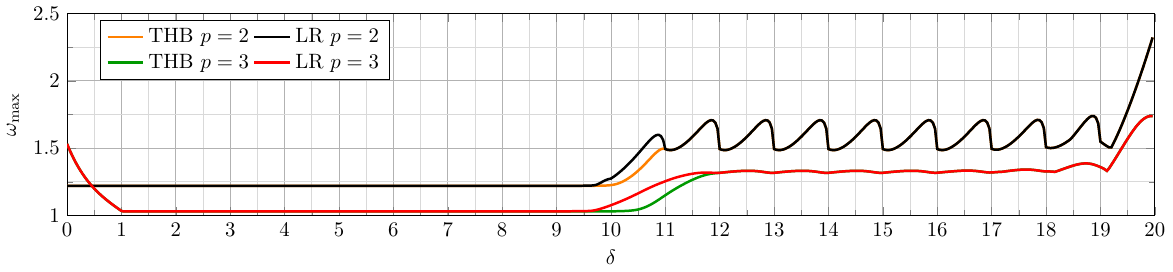}
    \caption{$\omega_{\text{max}}$ for THB- and LR-splines for $p=2$ and $p=3$.}
    \label{fig:omega_max_plot_conti}
\end{figure}

\rev{The continuous-trimming study in Fig.~\ref{fig:omega_max_plot_conti} provides further insight into the refinement-level requirement of BLCR. For $p=2$, the plateau from the onset of trimming shows that a one-level difference between the coarse boundary functions and the refined interior is sufficient in the investigated configuration, even while boundary basis functions remain active. For $p=3$, in contrast, the initial decrease of $\omega_{\max}$ before the plateau is reached confirms the stronger boundary contribution observed for the cubic case: a one-level difference alone is not sufficient while the most restrictive boundary functions remain active. Once a sufficiently large portion of the exterior boundary region has been trimmed away, $\omega_{\max}$ becomes independent of $\delta$ until the trimming boundary reaches the refined basis functions. In both cases, the plateau indicates that the coarse trimmed functions no longer govern the upper spectral limit. Instead, $\omega_{\max}$ is controlled by the refined untrimmed interior.}

\subsubsection{Spectral Bounds for Reduced-Support Basis Functions}
\label{subsub:Rayleigh}

\rev{The elevated eigenfrequencies associated with the exterior basis functions of open knot vectors were already identified by Adam et al.~\cite{ADAM2015581}. As discussed above, knot-exact trimming, unclamping, and arbitrary trimming can all modify the active support of basis functions, although the affected functions and the remaining support depend on the particular configuration. The purpose of the following analysis is therefore not to introduce a separate trimming mechanism, but to identify why reduced-support basis functions adjacent to a trimming boundary can remain critical after the exterior open-knot functions have been removed, and why refinement of these functions can again increase the upper spectral limit.}

To explain this behavior and clarify why BLCR alleviates it, we bound the largest eigenvalue using the Rayleigh quotient together with the modified Gershgorin estimate from Eq.~\ref{eq:mod_gersh}. Consider the generalized eigenvalue problem
\begin{equation}
    \mathbf{K}\mathbf{q} = \omega^2 \mathbf{M}\mathbf{q},
\end{equation}
where $\mathbf{K}$ and $\mathbf{M}$ denote the stiffness and lumped mass matrices. 
Assuming a symmetric stiffness matrix and a positive definite mass matrix, the Rayleigh quotient
\begin{equation}
    \rho(\boldsymbol{\phi}) =
    \frac{\boldsymbol{\phi}^T \mathbf{K} \boldsymbol{\phi}}
         {\boldsymbol{\phi}^T \mathbf{M} \boldsymbol{\phi}},
    \qquad \boldsymbol{\phi} \neq \mathbf{0},
\label{eq:Rayleigh_quotient}
\end{equation}
is bounded by the extremal eigenvalues
\begin{equation}
    \omega_{\min}^2 \le
    \rho(\boldsymbol{\phi})
    \le
    \omega_{\max}^2.
    \label{eq:min_max_bound}
\end{equation}
Choosing the canonical unit vector $\mathbf{e}_i$ for $\boldsymbol{\phi}$ in Eq.~\ref{eq:Rayleigh_quotient} yields the diagonal ratios 
\begin{equation}
    \rho(\boldsymbol{\mathbf{e}_i}) = Q_{ii}
    =
    \frac{\mathbf{e}_i^T \mathbf{K} \mathbf{e}_i}
         {\mathbf{e}_i^T \mathbf{M} \mathbf{e}_i}
    =
    \frac{K_{ii}}{M_{ii}},
\label{eq:RQ}
\end{equation}
which therefore provide computable lower bounds for the largest eigenvalue. Recalling the modified Gershgorin circle theorem from Eq.~\ref{eq:mod_gersh}, an upper bound for the largest eigenvalue is given by
\begin{equation}
    \omega_{\max}^2 
    \le 
    \left\| \mathbf{M}^{-1/2}\mathbf{K}\mathbf{M}^{-1/2} \right\|_\infty
    =
    \max_i 
    \left(
        \sum_{j=1}^{n}
        \frac{|K_{ij}|}{\sqrt{M_{ii}M_{jj}}}
    \right) = \bigl(\omega^{\mathrm{G,mod}}_{\max}\bigl)^2.
\end{equation}
Separating the diagonal contribution and using Eq.~\ref{eq:RQ} as lower bound yields
\begin{equation}
\max_i \frac{K_{ii}}{M_{ii}}
\le
\omega_{\max}^2
\le
\max_i
\left(
\frac{K_{ii}}{M_{ii}}
+
\sum_{j\neq i}
\frac{|K_{ij}|}{\sqrt{M_{ii}M_{jj}}}
\right).
\label{eq:bound}
\end{equation}
For brevity, we introduce the notation
\begin{equation}
    Q_{ij}=\frac{|K_{ij}|}{\sqrt{M_{ii}M_{jj}}},
\end{equation}
and we obtain
\begin{equation}
\max_i Q_{ii}
\le
\omega_{\max}^2
\le
\max_i
\left(
Q_{ii}
+
\sum_{j\neq i} Q_{ij}
\right).
\label{eq:bound2}
\end{equation}
To isolate the underlying mechanism, we consider the univariate setting. Since the bounds depend only on stiffness and lumped mass contributions and their scaling under refinement (see \ref{app:rayleigh_scaling}), the argument provides a consistent explanation also in higher dimensions. For basis functions $T_i(\xi)$, we get the stiffness and mass entries
\begin{equation}
K_{ij} =
\int_{\tilde{\Omega}_a}
T_i'(\xi)\,T_j'(\xi)\,\mathrm{d}\xi,
\qquad
M_{ii} =
\int_{\tilde{\Omega}_a}
T_i(\xi)\,\mathrm{d}\xi.
\end{equation}
The lumped mass is obtained by row-summing the consistent mass matrix. 
Owing to the partition-of-unity property of the basis, this is equivalent to integrating $T_i$ over its active support. 
The corresponding lower and upper frequency bounds associated with the $i$-th DOF are defined as
\begin{subequations}
\begin{align}
    \omega_i^{\mathrm{Q}}
    &=
    \sqrt{Q_{ii}},
    \label{eq:lower_upper_a}
    \\
    \omega_i^{\mathrm{G,mod}}
    &=
    \sqrt{
        Q_{ii}
        +
        \sum_{j\neq i} Q_{ij}
    }.
    \label{eq:lower_upper_b}
\end{align}
\label{eq:lower_upper}
\end{subequations}

\begin{figure}[H]
    \centering

    \begin{minipage}{0.9\textwidth}
        \centering
        \includegraphics[width=\linewidth]{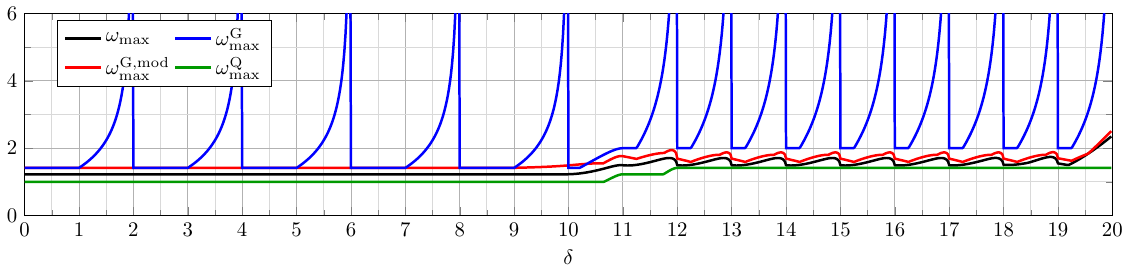}
        \subcaption{$p=2$}
        \label{fig:bounds_p2_total}
    \end{minipage}

    \vspace{3mm}

    \begin{minipage}{0.9\textwidth}
        \centering
        \includegraphics[width=\linewidth]{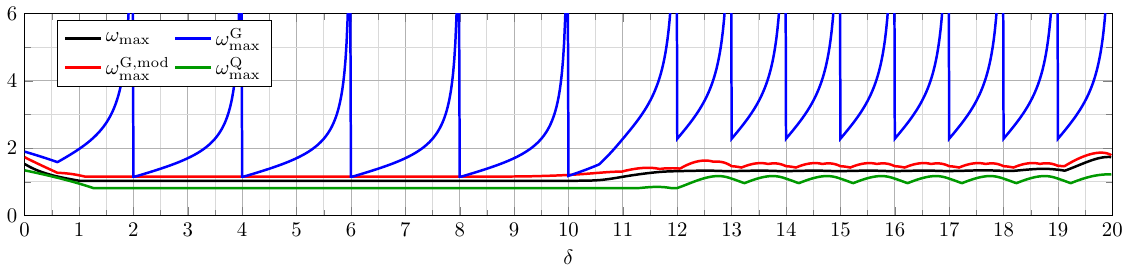}
        \subcaption{$p=3$}
        \label{fig:bounds_p3_total}
    \end{minipage}

    \caption{\rev{Maximum patch eigenfrequency and corresponding Rayleigh-based lower and Gershgorin upper bounds for symmetrically trimmed one-dimensional THB-spline discretizations. The trimming distance from both ends is denoted by $\delta$. The classical Gershgorin estimate is included for comparison but is not used in the subsequent BLCR argument.}}
    \label{fig:bounds_total}
\end{figure}

\rev{Fig.~\ref{fig:bounds_total} compares the maximum patch frequency with the global Rayleigh-based lower bound and the classical and modified Gershgorin upper bounds for $p=2$ and $p=3$. The classical Gershgorin estimate is retained for comparison. Its increasing overestimation for small active supports is well known, and it is therefore not used in the subsequent BLCR argument. The modified Gershgorin estimate remains substantially sharper over the relevant trimming range and is used for the DOF-wise comparison below.}

\rev{For $p=2$, $\omega_{\max}$ is independent of $\delta$ as long as no refined basis function is affected by trimming. For $p=3$, the corresponding plateau is reached for $\delta>1$. These plateaus confirm the interpretation obtained from Fig.~\ref{fig:omega_max_plot_conti}: within these trimming ranges, the coarse trimmed functions no longer govern the upper spectral limit, which is instead controlled by the refined untrimmed interior.}

\rev{To quantify this separation, let $\mathcal{T}_{\mathrm{int}}^{\mathrm{ref}}
\subset \mathcal{T}_{\mathrm{int}}$ denote the interior basis functions on the maximum refinement level. For the configurations considered here, a sufficient bound-separation condition for the coarse trimmed functions not to control the upper spectral limit is}
\begin{equation}
\rev{
    \underbrace{
    \max_{T_i\in\mathcal{T}_{\mathrm{cut}}}
    \omega_i^{\mathrm{G,mod}}
    }_{\displaystyle
    \omega_{\mathrm{cut,max}}^{\mathrm{G,mod}}
    }
    <
    \underbrace{
    \min_{T_j\in\mathcal{T}_{\mathrm{int}}^{\mathrm{ref}}}
    \omega_j^{\mathrm{Q}}
    }_{\displaystyle
    \omega_{\mathrm{int,min}}^{\mathrm{Q}}
    }.
}
\label{eq:blcr_bound_separation}
\end{equation}

\rev{If Eq.~\ref{eq:blcr_bound_separation} is satisfied, even the largest DOF-wise modified Gershgorin upper bound associated with a trimmed function remains below the smallest Rayleigh lower bound of the refined untrimmed interior functions. The two function classes are therefore separated in the bound analysis, and the coarse trimmed functions cannot constitute the limiting class for the maximum eigenfrequency.}

\rev{We evaluate this condition using representative maximally continuous B-spline basis functions at the coarse and refined levels. These are the two function classes central to the BLCR mechanism: coarse functions whose active support may be reduced by trimming and refined functions that remain in the interior. On a uniform mesh, maximally continuous interior B-splines are translations of one another, and varying $\delta$ across the support of one representative function covers the possible relative positions of the trimming boundary. Transition basis functions are deliberately excluded from this local comparison because their shapes depend on the LR- or THB-spline construction and would obscure the refinement-level mechanism common to both formulations. Their formulation-specific influence is captured directly by the preceding LR- and THB-spline investigations.}

\rev{Under uniform refinement, the lumped mass of an untrimmed interior basis function scales with $h$, whereas the corresponding stiffness entries scale with $h^{-1}$. Consequently, $Q_{ij}$ scales with $h^{-2}$ and the associated frequency bounds with $h^{-1}$. Element bisection therefore doubles the characteristic frequency bounds. The corresponding derivation, including the multidimensional setting and its limitations for trimmed functions, is given in~\ref{app:rayleigh_scaling}.}

\begin{figure}[H]
\centering

\begin{minipage}[t]{0.49\textwidth}
    \centering
    \includegraphics[width=\linewidth]{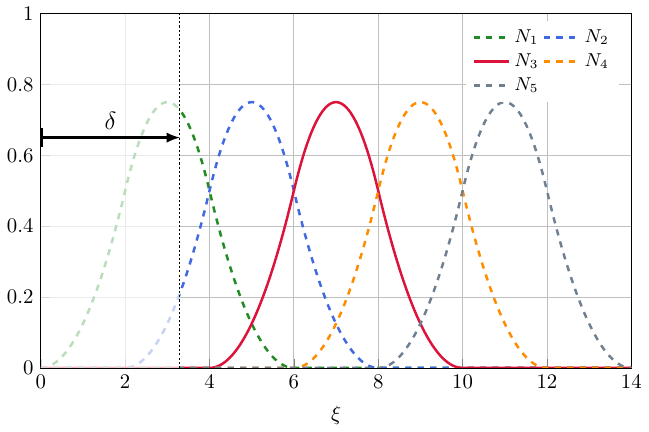}
    \subcaption{Basis functions, level~0}
    \label{fig:lvlcmp_p2_a}
\end{minipage}
\hfill
\begin{minipage}[t]{0.49\textwidth}
    \centering
    \includegraphics[width=\linewidth]{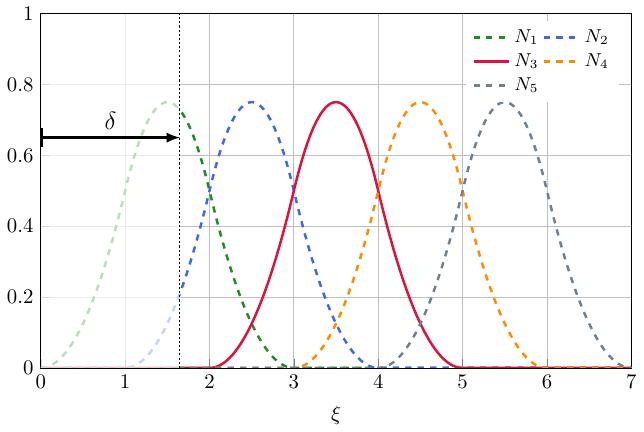}
    \subcaption{Basis functions, level~1}
    \label{fig:lvlcmp_p2_b}
\end{minipage}

\vspace{4mm}

\begin{minipage}[t]{0.49\textwidth}
    \centering
    \includegraphics[width=\linewidth]{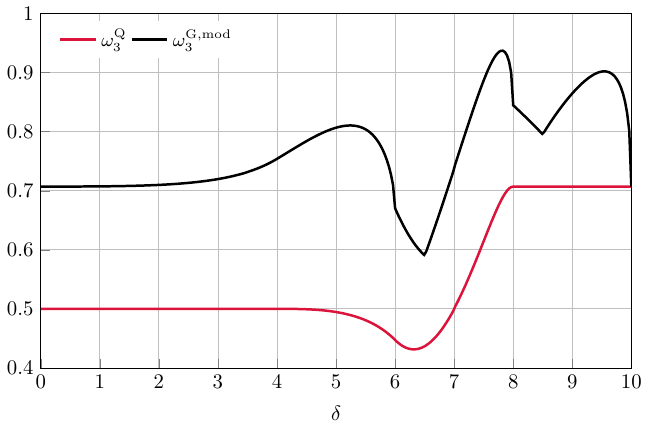}
    \subcaption{DOF-wise frequency bounds, level~0}
    \label{fig:lvlcmp_p2_bounds_l0}
\end{minipage}
\hfill
\begin{minipage}[t]{0.49\textwidth}
    \centering
    \includegraphics[width=\linewidth]{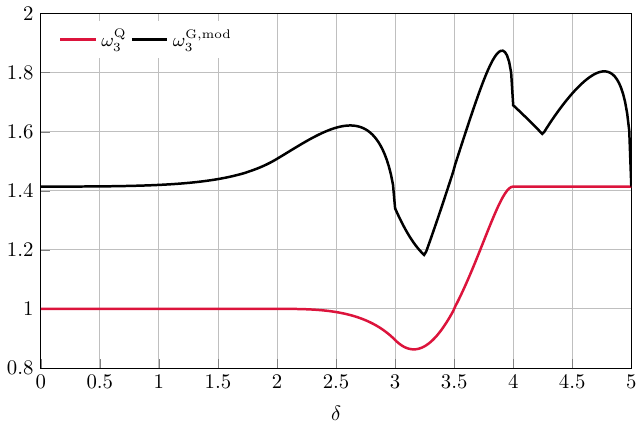}
    \subcaption{DOF-wise frequency bounds, level~1}
    \label{fig:lvlcmp_p2_bounds_l1}
\end{minipage}

\caption{\rev{Representative maximally continuous B-spline basis function and corresponding DOF-wise frequency bounds for $p=2$. Left: level~0 (coarse). Right: level~1 (refined). The red curve denotes the Rayleigh-based lower bound $\omega_3^{\mathrm{Q}}$, and the black curve the modified Gershgorin upper bound $\omega_3^{\mathrm{G,mod}}$.}}
\label{fig:level_compare_p2}
\end{figure}

\rev{For $p=2$, consider the representative coarse basis function shown in Fig.~\ref{fig:lvlcmp_p2_a}. In the untrimmed state ($\delta=0$), its contribution is bounded by $\omega_3^{\mathrm{Q}}=0.5$ and $\omega_3^{\mathrm{G,mod}}=\sqrt{2}/2\approx0.707$. As the trimming boundary moves through its support, both bounds vary with the remaining active support, as shown in Fig.~\ref{fig:lvlcmp_p2_bounds_l0}. In the knot-exact configuration relevant to the preceding investigations, at least one basis function adjacent to the trimming boundary has reduced support, and its Rayleigh lower bound reaches $\omega_3^{\mathrm{Q}}=\sqrt{2}/2\approx0.707$. Under uniform refinement, the corresponding frequency bounds scale by a factor of two, yielding $\omega_3^{\mathrm{Q}}=\sqrt{2}\approx1.414$ for the analogous reduced-support configuration in Fig.~\ref{fig:lvlcmp_p2_bounds_l1}.}

\rev{Under BLCR, the potentially trimmed functions remain on the coarse level. Over the relevant trimming range in Fig.~\ref{fig:lvlcmp_p2_bounds_l0}, their modified Gershgorin upper bound does not exceed approximately $0.94$. The refined untrimmed interior functions, represented by the $\delta=0$ state in Fig.~\ref{fig:lvlcmp_p2_bounds_l1}, have a Rayleigh lower bound of $1$. Hence, Eq.~\ref{eq:blcr_bound_separation} is satisfied. The coarse trimmed functions therefore lie entirely below the refined interior in the DOF-wise bound comparison. For the present $p=2$ configuration, $\omega_{\max}$ is confined to $[1,\sqrt{2}]$, compared with $[\sqrt{2},1.88]$ for the uniformly refined trimmed B-spline configuration.}

\begin{figure}[H]
\centering

\begin{minipage}[t]{0.49\textwidth}
    \centering
    \includegraphics[width=\linewidth]{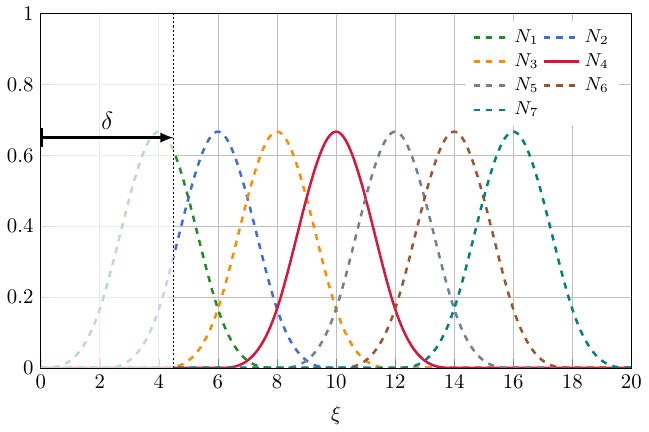}
    \subcaption{Basis functions, level~0}
    \label{fig:lvlcmp_p3_a}
\end{minipage}
\hfill
\begin{minipage}[t]{0.49\textwidth}
    \centering
    \includegraphics[width=\linewidth]{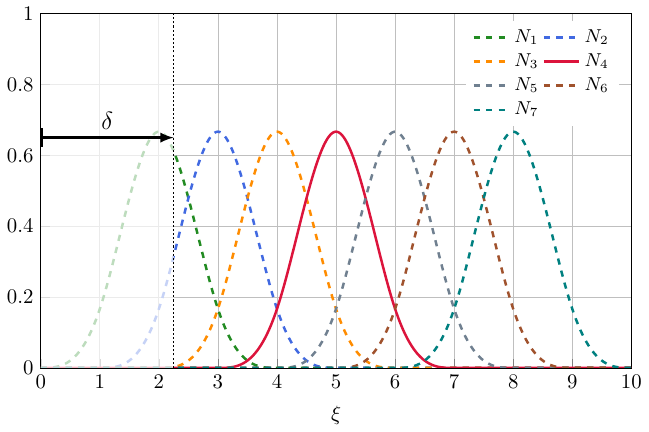}
    \subcaption{Basis functions, level~1}
    \label{fig:lvlcmp_p3_b}
\end{minipage}

\vspace{4mm}

\begin{minipage}[t]{0.49\textwidth}
    \centering
    \includegraphics[width=\linewidth]{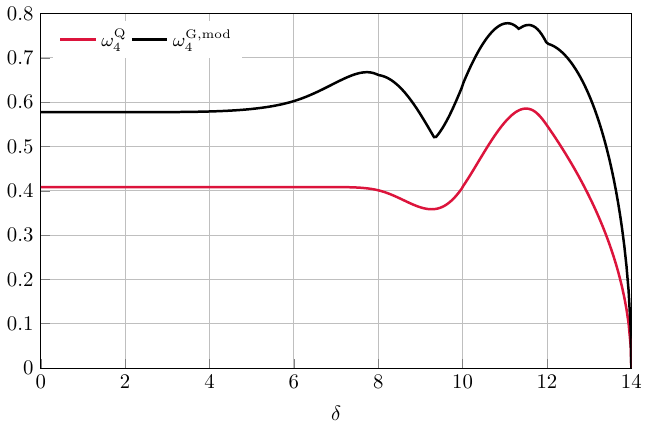}
    \subcaption{DOF-wise frequency bounds, level~0}
    \label{fig:lvlcmp_p3_bounds_l0}
\end{minipage}
\hfill
\begin{minipage}[t]{0.49\textwidth}
    \centering
    \includegraphics[width=\linewidth]{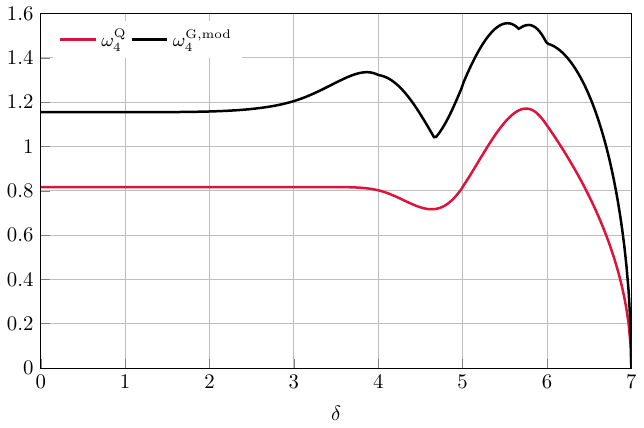}
    \subcaption{DOF-wise frequency bounds, level~1}
    \label{fig:lvlcmp_p3_bounds_l1}
\end{minipage}

\caption{\rev{Representative maximally continuous B-spline basis function and corresponding DOF-wise frequency bounds for $p=3$. Left: level~0 (coarse). Right: level~1 (refined). The red curve denotes the Rayleigh-based lower bound $\omega_4^{\mathrm{Q}}$, and the black curve the modified Gershgorin upper bound $\omega_4^{\mathrm{G,mod}}$.}}
\label{fig:level_compare_p3}
\end{figure}

\rev{For $p=3$, the same separation is observed. In the untrimmed state ($\delta=0$), the representative coarse basis function is bounded by $\omega_4^{\mathrm{Q}}=\sqrt{6}/6\approx0.408$ and $\omega_4^{\mathrm{G,mod}}=\sqrt{3}/3\approx0.577$. As its active support is reduced by trimming, both bounds increase. Over the relevant coarse trimming range, the modified Gershgorin upper bound remains below approximately $0.778$. The refined untrimmed interior functions, in contrast, have a Rayleigh lower bound of $\sqrt{6}/3\approx0.816$ and a modified Gershgorin upper bound of $2\sqrt{3}/3\approx1.155$. Hence, Eq.~\ref{eq:blcr_bound_separation} is again satisfied. For the investigated BLCR configuration, $\omega_{\max}$ is therefore confined to $[0.816,1.155]$, compared with $[0.963,1.556]$ for uniformly refined trimmed B-splines.}

\rev{The positive separation in Eq.~\ref{eq:blcr_bound_separation} explains why a one-level difference is sufficient for both $p=2$ and $p=3$ in the configurations investigated here. The size of this gap is configuration dependent. The present analysis therefore supports the observed one-level rule but does not establish a universal criterion for arbitrary geometries, polynomial degrees, or refinement patterns.}

\begin{remark}
\rev{Reduced active support may also affect the opposite end of the spectrum. For strongly trimmed basis functions, the diagonal Rayleigh quantity $Q_{ii}=K_{ii}/M_{ii}$ can become arbitrarily small, consistent with the low-frequency pathologies reported for trimmed and immersed row-sum-lumped discretizations \cite{Bioli2025,Voet_lumping_stabilization}.}
\end{remark}

\subsection{Investigations in 2D}
\label{2D_boundary_trimmed}
In this section we extend the refinement and trimming scenarios to the bivariate case. 
For all bivariate setups we use expanded open knot vector patches with at least \(p-1\) void element rows 
and columns along the boundary. Since trimming in two dimensions can become arbitrarily 
complex, we start with a simple square geometry and then consider a rotated variant. 
After reproducing the same studies as in the 1D case, we conclude with a complex example 
and demonstrate how BLCR can significantly reduce 
the maximum patch eigenfrequency and thus increase the critical time step in explicit dynamics compared to a globally refined patch.

\subsubsection{Symmetric Refinement from Center to Boundary for a Square}
\begin{figure}[H]
    \centering
    \renewcommand{\arraystretch}{1.4}

    \newcommand{\LeftBox}[1]{\makebox[0.5\textwidth][c]{#1}}
    \newcommand{\RightBox}[1]{\makebox[0.5\textwidth][c]{#1}}

    \begin{tabular}{@{} c c @{}}
        \LeftBox{
            \subfloat[][Refined Mesh at $w=16$.]{%
                \includegraphics[width=0.48\textwidth]{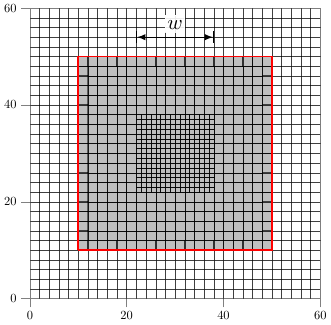}}
        } &
        \RightBox{
            \subfloat[][Refined Mesh at $w=40$.]{%
                \includegraphics[width=0.48\textwidth]{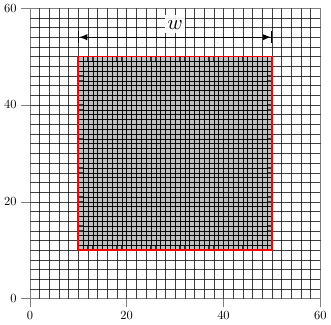}}
        }
    \end{tabular}
    \caption{Square with increasing refinement area.}
    \label{fig:square}
\end{figure}

\begin{figure}[H]
    \centering
    \includegraphics[width=\textwidth]{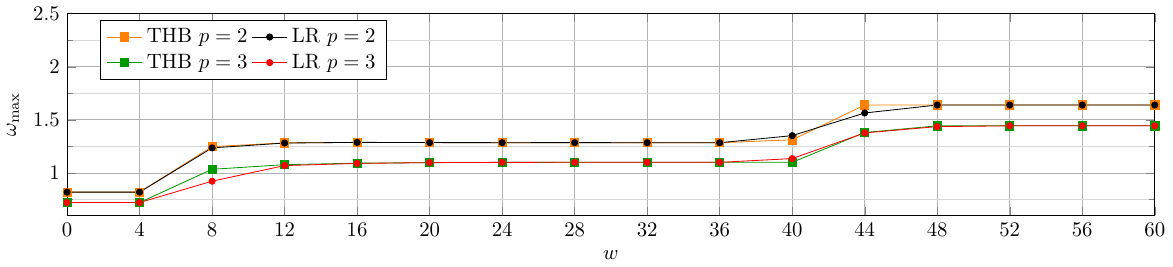}
    \caption{$\omega_{\text{max}}$ for THB- and LR-splines for $p=2$ and $p=3$.}
    \label{fig:square_omega}
\end{figure}

We start our investigation with a patch consisting of \(30 \times 30\) elements, each of size \(2 \times 2\). 
We trim off five element rows and columns, resulting in \(20 \times 20\) active elements. 
Fig.~\ref{fig:square} shows the setup. The center of the patch is refined symmetrically, controlled by the refinement 
width \(w\). Fig.~\ref{fig:square_omega} shows the maximum patch eigenfrequency for LR- and THB-splines for \(p=2\) and \(p=3\). 
A refinement width of \(w=4\) corresponds to an area of \(2 \times 2\) elements and does not fully contain 
a basis function. Consequently, the first effective refinement occurs for \(w \ge 8\). 

Similar to the 1D case, once local refinement is present, the maximum patch eigenfrequency becomes largely 
independent of the refinement width. For \(w=40\) (see Fig.~\ref{fig:square}b), all basis functions in the interior 
are refined. Increasing \(w\) further causes refinement of trimmed basis functions, \rev{which raises the upper spectral limit as the refined reduced-support functions become critical, analogous to the 1D observations.}

\subsubsection{Symmetric Refinement from Center to Boundary for a  Rotated Square}
\begin{figure}[H]
    \centering
    \renewcommand{\arraystretch}{1.4}

    \newcommand{\LeftBox}[1]{\makebox[0.5\textwidth][c]{#1}}
    \newcommand{\RightBox}[1]{\makebox[0.5\textwidth][c]{#1}}

    \begin{tabular}{@{} c c @{}}
        \LeftBox{
            \subfloat[][Refined Mesh at $w=16$.]{%
                \includegraphics[width=0.48\textwidth]{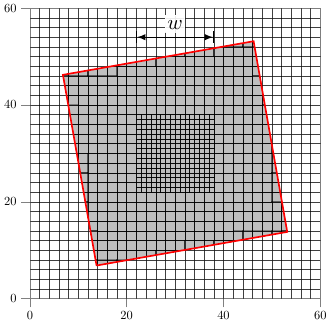}}
        } &
        \RightBox{
            \subfloat[][Refined Mesh at $w=36$, where the refined basis functions are trimmed.]{%
                \includegraphics[width=0.48\textwidth]{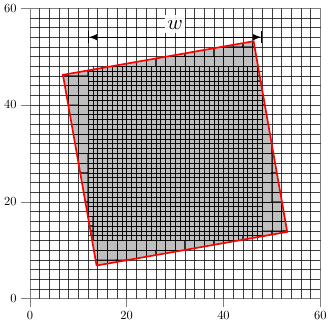}}
        }
    \end{tabular}

    \caption{Rotated square with increasing refinement area.}
    \label{fig:square_rot}
\end{figure}

\begin{figure}[H]
    \centering
    \includegraphics[width=\textwidth]{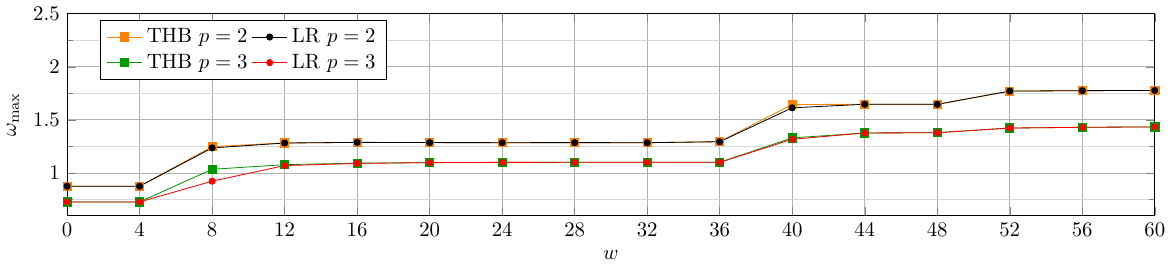}
    \caption{$\omega_{\text{max}}$ for THB- and LR-splines for $p=2$ and $p=3$.}
    \label{fig:square_rot_omega}
\end{figure}

We now use the previous setup and rotate the trimming curve by \(10^\circ\), as shown in Fig.~\ref{fig:square_rot}. 
The behavior of the maximum patch eigenfrequency as shown in Fig.~\ref{fig:square_rot_omega} is essentially the same as in the unrotated case. 
The only difference is that two distinct jumps appear at the upper end of the spectrum. As shown in 
Fig.~\ref{fig:square_rot}b, for \(w=36\) some refined basis functions are slightly trimmed. Increasing \(w\) causes 
the first jump, as refined basis functions along the square edge become trimmed. The second jump occurs for 
\(w=52\), where the refined corner basis functions of the rectangle are trimmed.

\subsubsection{Continuous Trimming from Boundary to Center for a Square}
\begin{figure}[H]
    \centering
    \includegraphics[width=0.48\textwidth]{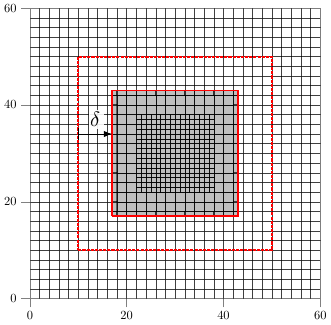}
    \caption{Square with fixed center refinement. Symmetric trimming distance is denoted by $\delta$.}
    \label{fig:square_trim}
\end{figure}

\begin{figure}[H]
    \centering
    \includegraphics[width=\textwidth]{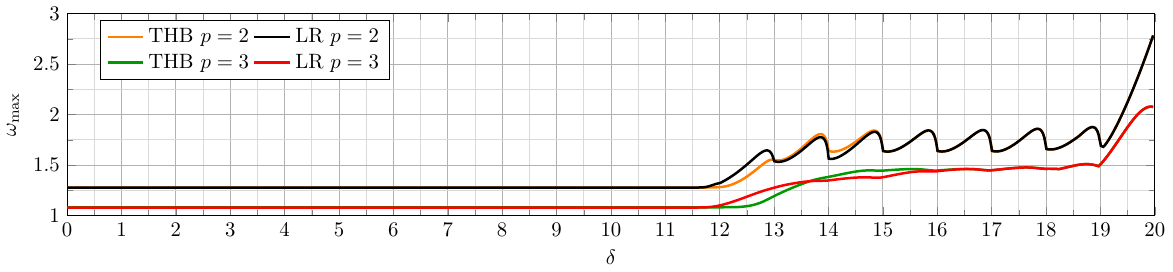}
    \caption{$\omega_{\text{max}}$ for THB- and LR-splines for $p=2$ and $p=3$.}
    \label{fig:square_trim_omega}
\end{figure}

We now consider an embedded rectangle of size \(40 \times 40\) with a fixed refinement width $w=16$ for LR- and THB-splines with \(p=2\) and \(p=3\), as shown in Fig.~\ref{fig:square_trim}. The rectangle is then trimmed 
symmetrically, with the trimming distance denoted by \(\delta\). Fig.~\ref{fig:square_trim_omega} shows the maximum patch 
eigenfrequency as a function of \(\delta\). The observed behavior is similar to the 1D case. As long as \(\delta < 12\), no refined basis functions 
are trimmed, and the maximum patch eigenfrequency remains independent of \(\delta\). Since we start 
with an expanded patch containing void boundary elements, no boundary effect is visible for \(p=3\). 
Once a refined basis function is trimmed, the maximum patch eigenfrequency increases, with the LR-splines showing a slightly earlier rise compared to THB-splines. \rev{As in 1D, the plateau before this transition spans continuously varying, non-knot-exact trimming configurations and demonstrates that coarse trimmed functions do not govern $\omega_{\max}$ while the refined region remains uncut.} For \(\delta > 19\), a strong increase in \( \omega_{\max} \) is observed, as fewer than \(2 \times 2\) refined elements and thus less than one full basis function support remain. This increase is consistent with the observations in the 1D case. While such small patches are less likely to occur, especially when analysis suitability is considered, this behavior aligns with our previous study \cite{hollweck_LR_THB_2026}, in particular regarding trimmed geometries with narrow ribs.

\subsubsection{Continuous Trimming from Boundary to Center for a Rotated Square}
\begin{figure}[H]
    \centering
    \includegraphics[width=0.48\textwidth]{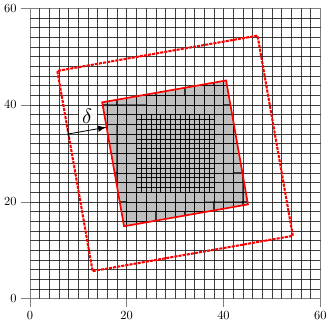}
    \caption{Rotated square with fixed center refinement. Symmetric trimming distance is denoted by $\delta$.}
    \label{fig:square_rot_trim}
\end{figure}

\begin{figure}[H]
    \centering
    \includegraphics[width=\textwidth]{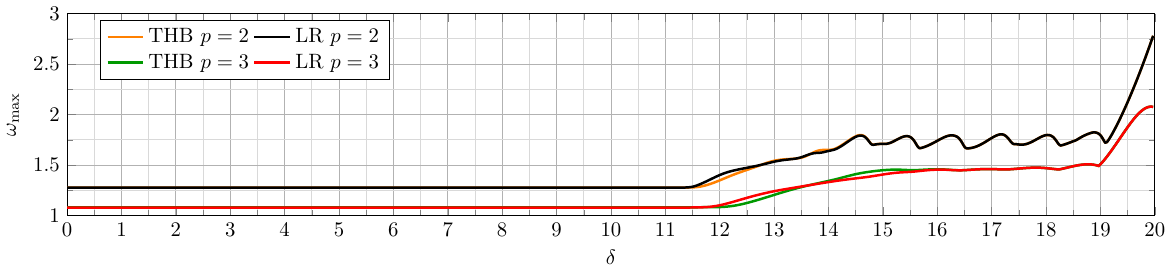}
    \caption{$\omega_{\text{max}}$ for THB- and LR-splines for $p=2$ and $p=3$.}
    \label{fig:square_rot_trim_omega}
\end{figure}

We now consider the same setup as in the previous section, with the only difference that the square is rotated by \(10^\circ\), as shown in Fig.~\ref{fig:square_rot_trim}. Fig.~\ref{fig:square_rot_trim_omega} shows the maximum patch eigenfrequency as a function of \(\delta\). The observed behavior is identical to the previous setup as long as no refined basis function is trimmed. Once a refined basis function is trimmed, the maximum patch eigenfrequency again increases, with LR-splines showing a slightly earlier rise compared to THB-splines. In contrast to the unrotated case, the transitions become smoother and the differences between LR- and THB-splines practically vanish. \rev{This confirms in a genuinely non-knot-exact two-dimensional setting that the relevant distinction for the upper spectral limit is whether reduced-support basis functions belong to the refined or the coarse level, rather than whether the trimming boundary coincides with element edges.}

\subsubsection{Complex Trimmed Geometry}
\begin{figure}[H]
    \centering

    \begin{tabular}{c c}
        \subfloat[][B-spline]{\includegraphics[width=0.48\textwidth]{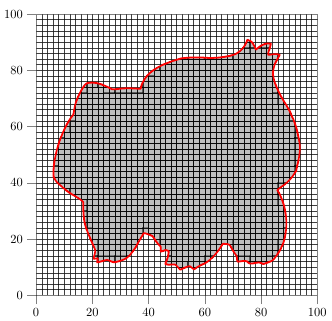}} &
        \subfloat[][globally refined B-spline]{\includegraphics[width=0.48\textwidth]{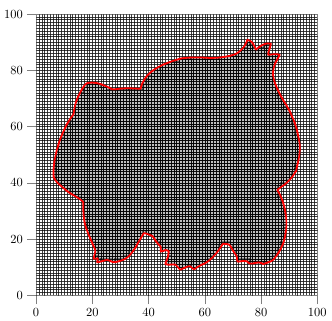}} \\[1em]
        \subfloat[][BLCR $p=2$]{\includegraphics[width=0.48\textwidth]{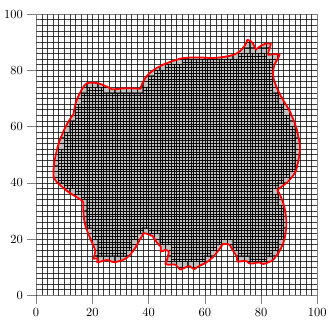}} &
        \subfloat[][BLCR $p=3$]{\includegraphics[width=0.48\textwidth]{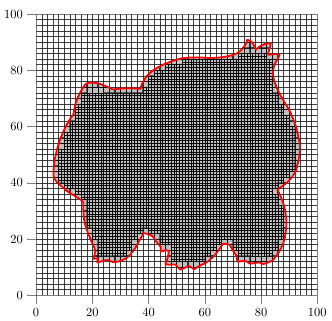}}
    \end{tabular}

    \caption{Complex geometry with different meshes.}
    \label{fig:pokemon}
\end{figure}

Fig.~\ref{fig:pokemon}a shows a complex trimmed geometry together with its globally refined counterpart in
Fig.~\ref{fig:pokemon}b. As already observed, global refinement also refines trimmed basis functions, \rev{thereby increasing the upper spectral limit associated with reduced-support functions.} To avoid this effect, we employ LR- and THB-splines and restrict refinement exclusively to active, untrimmed basis functions. This procedure results in identical meshes for LR- and THB-splines. The meshes, however, differ between polynomial degrees. With decreasing element size, these differences become negligible, as illustrated for \(p=2\) in Fig.~\ref{fig:pokemon}c and for \(p=3\) in Fig.~\ref{fig:pokemon}d. The only visible discrepancy is located in the vicinity of the coordinate \((20,70)\).
We now examine the maximum patch eigenfrequency for all configurations shown in Fig.~\ref{fig:pokemon}, using LR-
and THB-splines for the BLCR cases in Figs.~\ref{fig:pokemon}c and d. Between the globally refined configurations in Figs.~\ref{fig:pokemon}a and b, the maximum patch eigenfrequency increases - as expected - by approximately a factor of two, see Tab.~\ref{tab:pokemon}. In contrast, the BLCR configurations exhibit a substantial reduction in the maximum patch
eigenfrequency. This reduction corresponds to an increase of the critical time step of approximately
\(42\%\) for \(p=2\) and \(38\%\) for \(p=3\).

\begin{table}[H]
    \centering
    \caption{Maximum eigenfrequency $\omega_{\max}$ and resulting increase in critical time step for different refinement strategies.}
    \setlength{\tabcolsep}{14pt}
    \renewcommand{\arraystretch}{1.3}

    \begin{tabular}{|l||c|c|}
        \hline
        \textbf{} & \textbf{$p=2$} & \textbf{$p=3$} \\
        \hline\hline
        $\omega_{\max}$ for B-splines Level 0 & 0.87 & 0.75 \\ \hline
        $\omega_{\max}$ for B-splines Level 1 & 1.82 & 1.50 \\ \hline
        $\omega_{\max}$ for BLCR LR-splines & 1.28 & 1.09 \\ \hline
        $\omega_{\max}$ for BLCR THB-splines & 1.28 & 1.09 \\ \hline \hline
        $\Delta t_{\rm crit}^{\rm BLCR} / \Delta t_{\rm crit}^{\rm Level\,1}$ & 1.42 & 1.38 \\ \hline
    \end{tabular}
    \label{tab:pokemon}
\end{table}

\section{Industrial Application}
\label{sec:sec4}

In this section, we conduct a numerical study to validate our findings regarding the critical time step. In our analysis we use a cross bowl example from the field of sheet metal forming \cite{Wagner2026FEM}. This explicit simulation encompasses the complexities typically encountered in industrial applications, including contact, plasticity, large deformations, and high stress gradients. We use the isogeometric Reissner-Mindlin shell formulation implemented in LS-DYNA \cite{hallquist2006lsdyna}. For details on the formulation and implementation, see \cite{benson_isogeometric_2010, leidinger_explicit_2019, leidinger_Diss, DU2022109844, DU2024103728}. Furthermore, we use B\'ezier extraction as described in \cite{hollweck_LR_THB_2026} to incorporate the different spline formulations into LS-DYNA.

\subsection{Model setup for the Cross Bowl}

\begin{figure}[H]
  \centering
  \includegraphics[width=.99\textwidth]{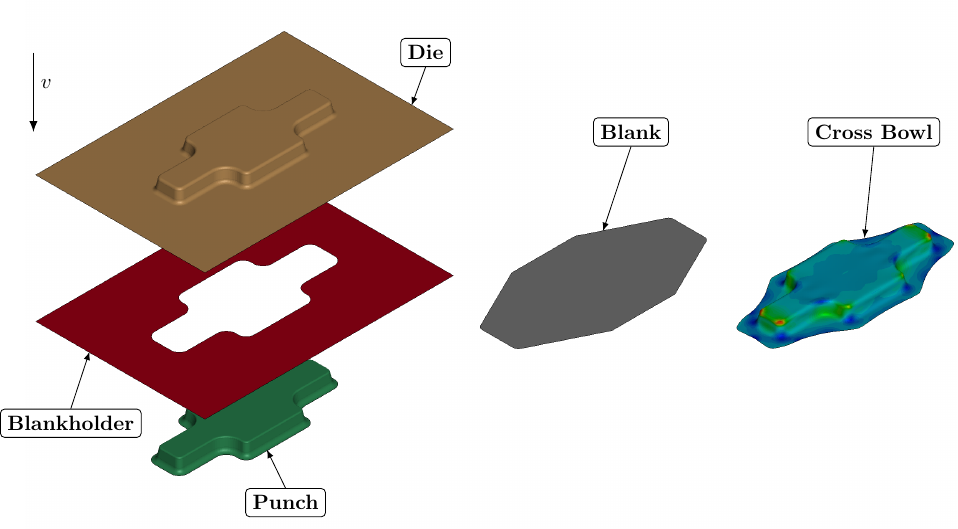}
  \caption{Cross Bowl forming setup including die, blankholder, punch, and blank. The blank is placed between the die and the blankholder. The die moves with prescribed velocity~\(v\) and pushes the blank into the die, forming the cross-bowl geometry.}
  \label{fig:SRail_tools}
\end{figure}

The cross bowl geometry consists of a curved sheet metal blank with a cross shape, which is obtained by drawing the blank into a die using a punch and blankholder, see Fig.~\ref{fig:SRail_tools}. The material is assumed to be isotropic and follows an elastoplastic behavior with isotropic hardening. The material properties are defined by a Young's modulus of 210\,GPa, a Poisson's ratio of $\nu$=0.3, and a yield stress of 275\,MPa. The hardening behavior is described by a measured curve relating true strain and true stress. The punch, die, and blankholder are modeled as rigid surfaces to reduce computational cost, a common approach in industrial applications to avoid expensive solid discretizations of the tools. A constant blank holder force of 130\,kN is applied to control material flow and prevent wrinkling during the forming process. The forming process is simulated using the explicit dynamic solver of LS-DYNA. The punch is fixed in all degrees of freedom, while the die is prescribed a vertical velocity curve to deform the blank into the die cavity. The blankholder is constrained to move only in the vertical direction, allowing it to apply the required clamping force. The blank is fully unconstrained. Coulomb friction is assumed at the contact interfaces between the blank and the tools, with a friction coefficient of $\mu=0.1$. A penalty-based, one-way surface-to-surface contact formulation is used, which is well suited for forming simulations \cite{hallquist2006lsdyna}.

\subsection{Preprocessing}

To apply BLCR in an industrial workflow, the initial geometry must be represented on a mesh that is one refinement level coarser than the target discretization. This implicitly requires that the coarser mesh, and thus the associated spline space, represents the geometry with sufficient accuracy. In practice, this is not a severe limitation for most industrial geometries. In particular, for completely flat geometries such as blanks in sheet metal forming, this requirement is trivially satisfied.

\begin{figure}
    \centering
    \renewcommand{\arraystretch}{1.4}

    \newcommand{\LeftBox}[1]{\makebox[0.5\textwidth][c]{#1}}
    \newcommand{\RightBox}[1]{\makebox[0.5\textwidth][c]{#1}}

    \begin{tabular}{@{} c c @{}}
        \LeftBox{
            \subfloat[][globally refined B-splines]{%
                \includegraphics[width=0.48\textwidth]{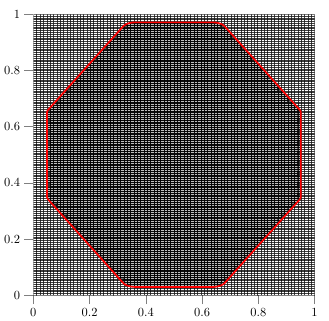}}
        } &
        \RightBox{
            \subfloat[][BLCR for LR- and THB-splines]{%
                \includegraphics[width=0.48\textwidth]{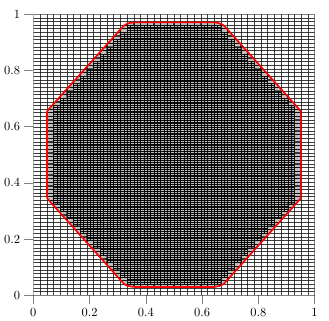}}
        }
    \end{tabular}

    \caption{Geometry for the blank with different meshes. For this simple geometry the locally refined meshes are the same for $p=2$ and $p=3$.}
    \label{fig:blank_ibc}
\end{figure}

Starting from the coarse representation, the BLCR strategy is applied by locally refining only those basis functions that are interior and untrimmed, while deliberately leaving the trimmed boundary basis functions on the coarse level. Figs.~\ref{fig:blank_ibc}a and b show the different meshes for global refinement and BLCR, respectively.

\subsection{Results}

We now demonstrate the effectiveness of the BLCR strategy for LR- and THB-splines with $p=2$ and $p=3$.
We start from a standard tensor-product mesh using B-splines with an initial element size of
$2.5\,\mathrm{mm}$. Subsequently, global refinement is applied for B-splines, whereas for LR- and
THB-splines the BLCR strategy is used, resulting in a target element size of $1.25\,\mathrm{mm}$.

We apply the isogeometric Reissner-Mindlin shell formulation implemented in LS-DYNA \cite{hallquist2006lsdyna, benson_isogeometric_2010}. As is standard for Reissner-Mindlin shells in explicit dynamics, the rotational inertia terms are scaled such that the rotational modes do not restrict the time step. Consequently, the critical time step is governed by the largest membrane eigenfrequency, allowing the results from the previous membrane setups to be directly transferred to the shell discretization.

The theoretical critical time steps for the initial configuration are summarized in Tab.~\ref{tab:dtcrit_crossbowl}.

\begin{table}[H]
    \centering
    \caption{Critical time step $\Delta t_{\text{crit}}$ for the cross-bowl benchmark and relative increase obtained with BLCR compared to the globally refined trimmed B-spline discretization.}
    \setlength{\tabcolsep}{14pt}
    \renewcommand{\arraystretch}{1.3}

    \begin{tabular}{|l||c|c|}
        \hline
        \textbf{} & \textbf{$p=2$} & \textbf{$p=3$} \\
        \hline\hline
        $\Delta t_{\text{crit}}$ for B-splines &
        $5.68 \times 10^{-7}$ &
        $6.85 \times 10^{-7}$ \\ \hline
        $\Delta t_{\text{crit}}$ for BLCR (LR/THB) &
        $7.50 \times 10^{-7}$ &
        $8.80 \times 10^{-7}$ \\ \hline \hline
        $\Delta t_{\text{crit}}^{\text{BLCR}} / \Delta t_{\text{crit}}^{\text{B-spline}}$ &
        1.32 &
        1.28 \\ \hline
    \end{tabular}
    \label{tab:dtcrit_crossbowl}
\end{table}

BLCR increases the admissible time step by approximately $32\%$ for $p=2$ and $28\%$ for $p=3$ compared to the globally refined trimmed B-spline discretization. The increase is lower than that observed in the previous section. \rev{The trimming configuration of the cross-bowl geometry is relatively simple and does not involve sharp corners or strongly irregular trimming features. The achievable time-step increase therefore remains configuration dependent.}
In industrial explicit simulations, the critical time step is often estimated and updated during the
analysis. In practice, this estimate is typically conservative and does not account for the removal
of the trimming-induced boundary effect achieved by BLCR. In the present study, we therefore keep the time step constant throughout the simulation. As the mesh deforms, the maximum patch eigenfrequency naturally changes. To account for this, a safety factor of $f_s = 0.8$ is applied to the selected time step from Tab.~\ref{tab:dtcrit_crossbowl}.

\rev{To quantify the agreement between the different setups, we define the maximum total-energy deviation normalized by the final reference energy over the time interval $[0,T]$ as}
\[
\varepsilon_E
=
\left\|
\frac{E_{\mathrm{BLCR}}(t) - E_{\mathrm{ref}}(t)}
     {E_{\mathrm{ref}}(T)}
\right\|_{L^\infty(0,T)},
\]
\rev{where $E_{\mathrm{ref}}(t)$ denotes the total energy of the globally refined B-spline reference solution and $T$ is the final simulation time. The constant normalization by $E_{\mathrm{ref}}(T)$ is used deliberately to avoid division by zero or very small reference energies at the beginning of the simulation.}

\begin{figure}[H]
    \centering
    \includegraphics[width=0.7\textwidth]{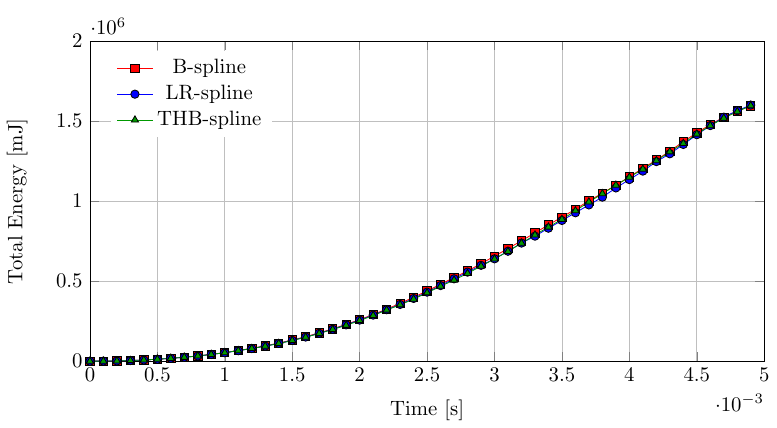}
    \caption{Evolution of the total energy for $p=2$.
    LR- and THB-splines are refined using the BLCR strategy, whereas B-splines are uniformly refined over the entire patch.}
    \label{fig:energy_plot_p2}
\end{figure}

\begin{figure}[H]
    \centering
    \includegraphics[width=0.7\textwidth]{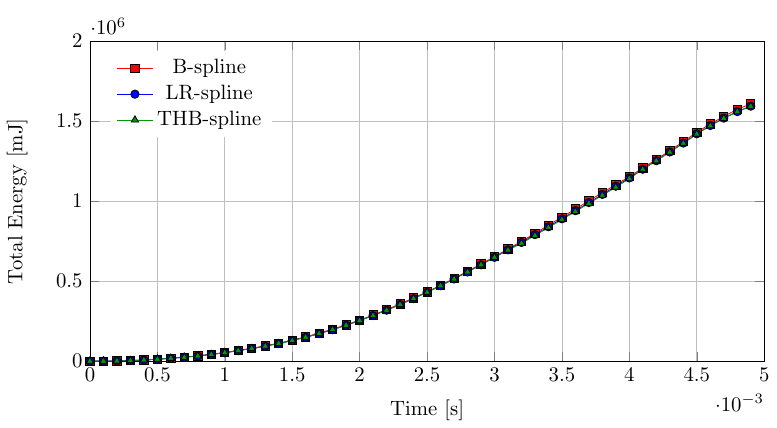}
    \caption{Evolution of the total energy for $p=3$.
    LR- and THB-splines are refined using the BLCR strategy, whereas B-splines are uniformly refined over the entire patch.}
    \label{fig:energy_plot_p3}
\end{figure}

\begin{figure}[H]
  \centering
  \includegraphics[width=0.65\textwidth]{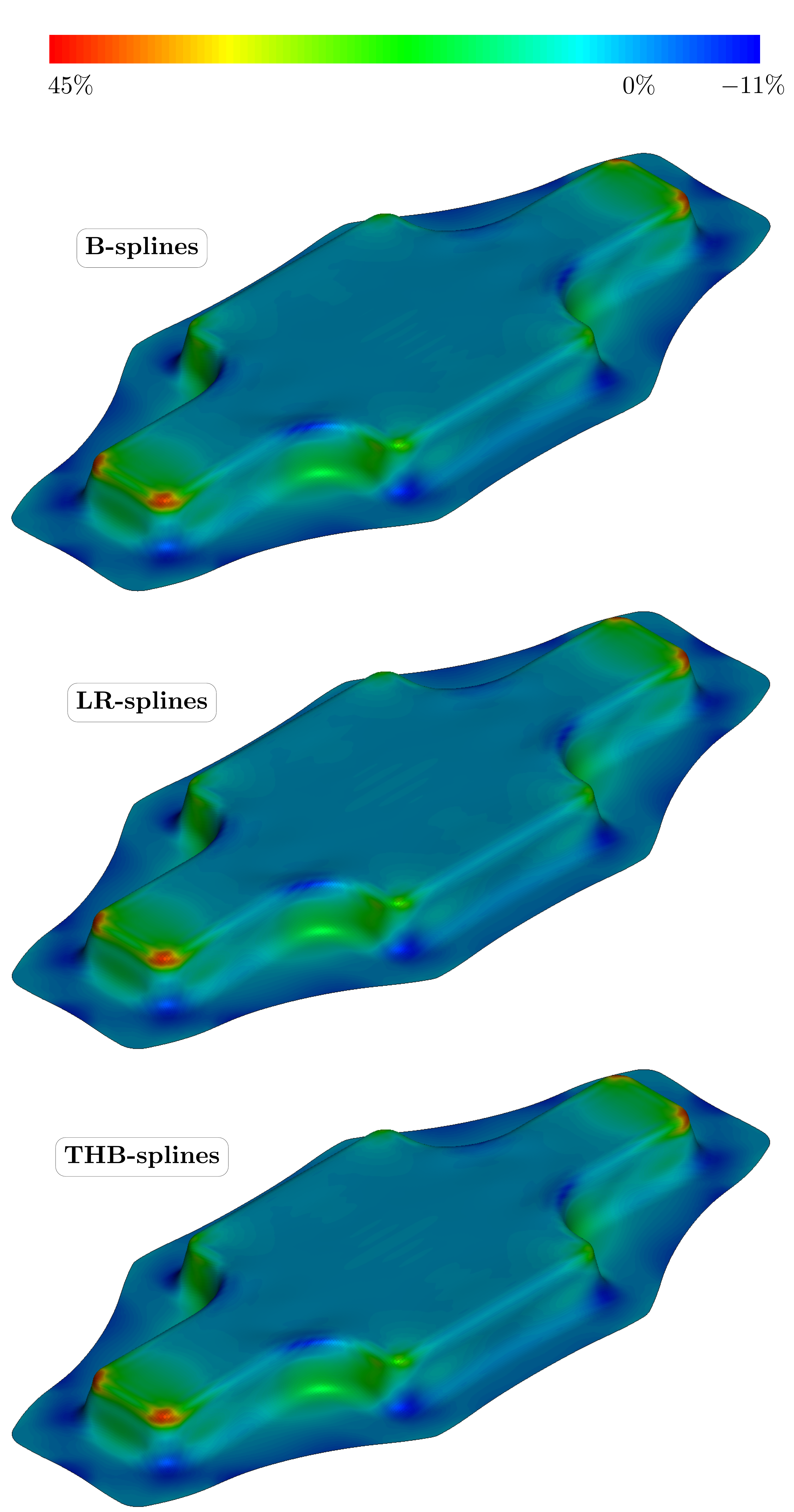}
  \caption{Change in shell thickness depicted in \%, where positive values (red) indicate thinning and negative values (blue) indicate thickening.}
  \label{fig:crossbowl_B_LR_THB}
\end{figure}

For $p=2$, the maximum deviation amounts to $\varepsilon_E = 1.1\%$ for THB-splines and $1.6\%$ for LR-splines. For $p=3$, it reduces to $1.0\%$ (THB) and $0.9\%$ (LR). The corresponding total-energy evolutions are shown in
Fig.~\ref{fig:energy_plot_p2} and Fig.~\ref{fig:energy_plot_p3}. No systematic drift or artificial energy accumulation is observed over the simulation time. Fig.~\ref{fig:crossbowl_B_LR_THB} shows the thickness distribution of the blank for B-, LR-, and THB-splines with $p=3$. No significant differences between the three simulations are observed.
\rev{These results demonstrate that, for the present forming benchmark, the increased critical time step obtained with BLCR is accompanied by close agreement with the globally refined reference solution in the investigated total-energy and thickness results.}

\subsection{Applicability and limitations}

The BLCR strategy deliberately maintains a coarser refinement level for basis functions adjacent to the trimming boundary. Consequently, the local approximation quality near trimmed edges is reduced compared to uniformly refined trimmed B-splines. The strategy therefore introduces a controlled reduction of boundary resolution in exchange for a smaller maximum eigenfrequency and a larger admissible time step.

In sheet metal forming applications, this limitation is often non-critical. The boundary of the blank typically represents scrap material that is removed after forming and does not belong to the final component. A reduced discretization accuracy in these regions may therefore have only a limited influence on the quantities of interest associated with the final part. In such settings, BLCR can improve the computational efficiency of locally refined discretizations. \rev{Adam et al. \cite{ADAM2015581} likewise identified metal-stamping configurations with deformation located away from the boundary as particular cases in which a coarser boundary resolution may be acceptable.}

In industrial practice, element sizes are typically selected as coarse as possible while satisfying accuracy requirements. In explicit dynamics, uniform global refinement approximately halves the critical time step due to the scaling of the maximum eigenfrequency, while simultaneously increasing the number of elements by a factor of four in two dimensions. This combination significantly raises the overall computational cost.

BLCR enables local refinement while avoiding the full time step penalty associated with refining trimmed basis functions.
In the investigated configurations, the reduction of the admissible time step after one refinement step is approximately 31\% (for $p=2$ and $p=3$), compared to about 50\% under uniform global refinement.
\rev{As a result, the time-step penalty associated with refinement is reduced.}

When combined with an adaptive refinement strategy, BLCR can further reduce the number of elements required to resolve localized phenomena, while limiting the impact on the critical time step. This provides additional potential for improving overall efficiency in explicit simulations.

For more complex multi-patch models, such as those encountered in crash simulations, the applicability of BLCR depends on the modeling objective. If accurate stress, strain, or damage prediction is required directly at trimmed boundaries, a coarser boundary discretization may be insufficient. In such cases, boundary refinement can be combined with local mass scaling of trimmed and refined basis functions in order to control the associated spectral outliers. This approach restores boundary resolution but introduces additional artificial mass, which may influence the dynamic response and therefore requires careful assessment.

\rev{In addition to the reduced boundary resolution, two distinct spectral-accuracy aspects must be considered. First, row-sum mass lumping itself modifies the discrete spectrum relative to the consistent-mass formulation. Its influence on the present forming benchmark may be comparatively limited because the simulation is conducted in a quasi-static regime, with kinetic energy remaining small relative to the internal energy. Since the mass matrix enters the equations of motion through the inertial term, a response dominated by plastic deformation, contact forces, and prescribed tool motion is expected to be less sensitive to the mass approximation than a wave propagation, vibration, or resonance problem.}

\rev{Second, trimming and the BLCR may alter individual frequencies in the lower part of the spectrum, as discussed in Sec.~\ref{sec:sec3}.} \rev{Such a lower-spectrum deviation alone does not imply a significant error in the transient response. Its relevance depends on whether the associated response is sufficiently excited and contributes to the quantities of interest. The strong contact-driven kinematic guidance in the present forming process may reduce the excitation of localized spectral components.} \rev{Moreover, strongly trimmed basis functions possess only a small active support, such that large motions of the associated control points may remain spatially localized and have only a limited influence on the global deformation field. Related behavior of so-called light control points has been reported for trimmed explicit isogeometric discretizations \cite{leidinger_Diss}. Nevertheless, excessive local velocities or displacements may still deteriorate the numerical solution or even lead to solver failure.}

\rev{The close agreement observed here therefore demonstrates that no detrimental influence is apparent in the investigated response quantities, but it does not exclude the occurrence of the previously discussed spectral-accuracy problems in other applications.} 

\rev{For applications in which accurate modal properties, wave propagation, or vibration response are of primary importance, consistent-mass formulations combined with suitable outlier-removal techniques provide an alternative where the associated computational cost is acceptable \cite{Voet_outlier}. Dedicated eigenvalue-stabilization techniques provide another route for suppressing trimming-induced high-frequency outliers and recovering feasible critical time-step sizes in immersed explicit dynamics \cite{Eisentraeger2024EVS,Burchner2026GEVS}.}

\rev{If a diagonal mass matrix is required for efficient explicit time integration, stabilization prior to mass lumping can substantially improve the low-frequency accuracy of strongly trimmed discretizations \cite{Voet_lumping_stabilization,Guarino2025}.} \rev{The polynomial-extension approach of Voet et al.\ specifically targets small cut elements and basis functions whose active support is confined to such elements. By modifying the discrete formulation prior to mass lumping, trimming-induced high-frequency modes of the consistent-mass formulation are prevented from appearing as spurious modes in the low-frequency spectrum after lumping. In the investigated knot-exact configurations, however, no partially cut elements are present and the corresponding polynomial-extension criterion does not become active, while noticeable lower-spectrum deviations are nevertheless observed. The present results therefore indicate that trimming-related lower-spectrum inaccuracies under row-sum lumping are not restricted to the small-cut configurations addressed by this stabilization. A further limitation of the extension procedure is that positivity of the spline basis is not generally preserved, such that the stabilized consistent mass matrix may contain negative entries and standard row-sum lumping is no longer guaranteed to produce a positive-definite diagonal mass matrix \cite{Voet_lumping_stabilization}.}

\rev{A related but conceptually different treatment of strongly trimmed degrees of freedom was proposed by Leidinger \cite{leidinger_Diss}. So-called light control points are identified from their small lumped masses and stabilized relative to neighboring stable control points. Although light-control-point instabilities and trimming-induced spurious low-frequency modes are both associated with strongly reduced active supports and small lumped masses, an equivalence between the two phenomena has not been established.}

\subsection{Cross-talk effects}

In sheet metal forming applications, small trimming features are rarely encountered. However, in other types of analyses, such as crashworthiness simulations, patches frequently contain narrow slits, holes, or cutouts. For immersed methods, these features may lead to the so-called cross-talk effect, first reported in \cite{coradello_hierarchically_2020} and recently studied in detail in \cite{Lian2025CPD}.

Small trimming features may separate the support of basis functions into several disconnected regions, which can introduce non-physical coupling between spatially separated parts of the domain and thereby reduce the predictive quality of the simulation. Due to their higher continuity and larger support, spline-based discretizations in isogeometric analysis are particularly sensitive to these effects \cite{Lian2025CPD}.

A complete separation of the support domain is referred to as \emph{Type~1 cross-talk}. In two and three dimensions, even partially trimmed supports may induce spurious traction forces, commonly classified as \emph{Type~2 cross-talk}. Local refinement of the trimming boundary has been proposed to eliminate Type~1 cross-talk \cite{coradello_hierarchically_2020}, while Type~2 cross-talk vanishes in the limit of successive $h$-refinement.

At first glance, this strategy appears contrary to the philosophy of the BLCR strategy. Nevertheless, since cross-talk effects are strictly localized near the trimming boundary, one may deliberately refine these regions and apply mass scaling to the affected basis functions in order to control the resulting spectral outliers. Alternatively, the BLCR strategy can be combined with the control point duplication approach proposed in \cite{Lian2025CPD} to mitigate Type~1 cross-talk.

\newpage

\newpage
\section{Conclusion and Outlook}
\label{sec:sec5}

\rev{While our previous research \cite{hollweck_LR_THB_2026} provided a general assessment of critical time steps for locally refined splines, the present study investigates the remaining upper-spectral limitation after the particularly restrictive open-knot boundary elements have been removed. The results show that trimmed basis functions with reduced active support adjacent to a trimming boundary can still increase the maximum eigenfrequency, even for knot-exact trimming without arbitrarily small cut cells. The relevant mechanism is therefore associated with the non-proportional changes in stiffness and row-sum-lumped mass contributions caused by support reduction.}

\rev{To control this effect, this study introduces the Boundary-Level-Constrained Refinement (BLCR) strategy. BLCR constrains the refinement level of basis functions whose support intersects the trimming boundary relative to that of the refined interior. In all configurations investigated in this work, maintaining these functions on a sufficiently coarse level prevents trimmed basis functions from governing the upper spectral limit, such that the maximum eigenfrequency is instead controlled by the refined untrimmed interior.}

To support these observations, we utilized Rayleigh quotients and a modified version of the Gershgorin circle theorem and assessed the resulting behavior through 1D and 2D studies as well as an industrial sheet metal forming benchmark. The main findings are summarized as follows:

\begin{itemize}

\item \textbf{Reduced-support contribution after boundary treatment:}
\rev{Even after the open-knot boundary elements have been removed, trimmed basis functions adjacent to a trimming boundary can still contribute significantly to the maximum eigenfrequency. An arbitrarily small cut fraction is not required to affect the upper spectral limit.}

\item \textbf{Penalty of boundary refinement:}
\rev{In the investigated 1D and 2D configurations, refining basis functions whose support intersects the trimming boundary produces a pronounced increase in $\omega_{\text{max}}$. Continuous trimming studies further show that $\omega_{\text{max}}$ remains constant while only coarse boundary functions are trimmed and increases once the trimming boundary reaches refined basis functions.}

\item \textbf{Frequency bounds:}
The maximum Rayleigh frequency $\omega^{\text{Q}}_{\text{max}}$ provides a rigorous lower bound, whereas the modified Gershgorin estimate $\omega_{\text{max}}^{\mathrm{G,mod}}$, derived from the $\infty$-matrix norm of the symmetrically scaled stiffness matrix, provides an upper bound for $\omega_{\text{max}}$.
\rev{For the investigated $p=2$ and $p=3$ configurations, the DOF-wise modified Gershgorin upper bound of the coarse trimmed functions lies below the Rayleigh lower bound of the refined untrimmed interior functions. This separation explains why a one-level difference is sufficient in these cases. Its magnitude is configuration dependent and does not establish a universal one-level criterion.}

\item \textbf{Effectiveness of BLCR:}
\rev{BLCR shifts the upper spectral limit from trimmed basis functions to the refined untrimmed interior. Compared to uniformly refined trimmed patches with the same interior resolution, the investigated BLCR configurations yield a substantially lower maximum eigenfrequency and therefore a larger critical time step. This benefit is accompanied by changes in individual lower frequencies, reflecting an application-dependent trade-off between time-step efficiency and spectral accuracy.}

\item \textbf{LR- vs. THB-splines:}
Both spline formulations are suitable for BLCR. Minor differences arise from their different refinement mechanisms. In one dimension, THB-splines preserve coarse-level basis-function shapes more strictly at refinement transitions, whereas LR-splines modify them through the splitting procedure. In the investigated two-dimensional configurations, these differences become negligible, and the resulting critical time steps are effectively identical for both formulations under BLCR.

\item \textbf{Industrial Application:}
\rev{The cross-bowl forming simulation conducted with LS-DYNA demonstrates that BLCR can provide a significant increase in the critical time step while maintaining close agreement with the globally refined reference solution in the investigated total-energy and thickness results. This supports the applicability of BLCR to the considered class of sheet metal forming problems, but does not constitute a general statement regarding spectral accuracy in nonlinear dynamics.}

\end{itemize}

\rev{Overall, the results demonstrate that BLCR can improve the time-step efficiency of trimmed explicit IGA when local refinement is required in the interior and the trimming boundary is not a primary region of approximation interest. The strategy preserves the row-sum-lumped diagonal mass matrix and does not modify the total mass. At the same time, maintaining a coarser boundary resolution introduces an application-dependent trade-off that must be considered when accurate response quantities are required near trimmed boundaries.}

\rev{BLCR specifically targets the upper spectral limit governing the critical time step and should not be interpreted as a general spectral stabilization technique. As shown in the spectral investigations, the refinement constraint can alter individual frequencies in the lower part of the spectrum. The relevance of these changes depends on the excitation and the quantities of interest and remains an open question, particularly for strongly nonlinear dynamics.}

Looking beyond BLCR, strategies for reducing upper-spectrum outliers in large-scale explicit solvers should respect essential practical constraints in order to avoid unnecessary computational overhead. In particular, the diagonal structure of the lumped mass matrix should be preserved, as it enables trivial inversion, minimal memory consumption, and component-wise acceleration updates without solving a coupled linear system, thereby supporting efficient parallel time integration.

Furthermore, practical approaches should avoid requiring global assembly of mass or stiffness matrices. Explicit solvers typically operate on element-level force evaluations and assemble only internal and external force vectors. In nonlinear simulations involving contact, plasticity, and large deformations, global stiffness matrices are usually neither assembled nor stored. Internal forces are evaluated directly from the current stress state, and element stiffness matrices are often not formed explicitly. Approaches that rely on such quantities therefore introduce additional complexity and computational overhead.

As an alternative to modifying the spline space as in the BLCR strategy, one may employ local mass scaling applied to trimmed basis functions. In this approach, stiffness entries remain unchanged while the masses are increased to reduce the corresponding eigenfrequencies. Although this may mitigate spectral outliers, it alters the total mass and inertia of the system. Determining suitable scaling factors and assessing their impact on highly nonlinear simulations, including contact, plasticity, and damage, remains an open problem.

In this context, eigenvalue bounds such as the Rayleigh quotient or Gershgorin circles may provide simple and robust heuristics for selecting appropriate scaling factors for trimmed basis functions. Such a mass-scaling approach may be preferable when local refinement is unavailable or not desired in practice.

\appendix
\section{Asymptotic Behavior of the Maximum Eigenfrequency under Mesh Refinement}
\label{app:rayleigh_scaling}
\rev{This appendix summarizes standard and elementary scaling arguments for the maximum discrete eigenfrequency $\omega_{\max}$ of the generalized eigenvalue problem $\mathbf{K}\mathbf{q}=\lambda\,\mathbf{M}\mathbf{q}$ under uniform mesh refinement. The purpose is not to derive a new estimate, but to make explicit the relations used in the main text for the present B-spline setting.}
\rev{For elliptic operators involving first-order spatial derivatives (e.g., the Laplace operator), representative stiffness and row-sum-lumped mass contributions associated with interior B-spline basis functions scale as}
\[
K_{ii} \sim h^{d-2},
\qquad
M_{ii} \sim h^{d},
\]
\rev{where \(d\) denotes the spatial dimension. Their ratio therefore scales as $h^{-2}$, corresponding to a characteristic frequency scaling of $h^{-1}$. In the following, this scaling is first derived for the DOF-wise Rayleigh quantities and subsequently related to the maximum eigenfrequency of the full discrete system.}

\subsection{Bounds on the Maximum Eigenfrequency}
The maximum discrete eigenvalue $\lambda_{\max} = \omega_{\max}^2$ of the generalized eigenvalue problem $\mathbf{K}\mathbf{q} = \lambda \mathbf{M}\mathbf{q}$ is defined by the supremum of the Rayleigh quotient. A lower bound is provided by the maximum of the diagonal ratios of the stiffness and mass matrices:
\begin{equation}
\max_i \frac{K_{ii}}{M_{ii}} := \max_i Q_{ii} \le
\omega_{\max}^2 = \max_{\mathbf{q} \neq \mathbf{0}}
\frac{\mathbf{q}^T \mathbf{K} \mathbf{q}}
     {\mathbf{q}^T \mathbf{M} \mathbf{q}} .
\end{equation}
Defining the Rayleigh frequency as $\omega_{i}^{\mathrm{Q}} = \sqrt{Q_{ii}}$, the lower bound can be expressed in terms of these frequencies. For a discretization using B-spline basis functions and a lumped mass formulation, the diagonal entries correspond to the integrals of the squared derivatives and the basis functions themselves. In the following, we analyze how these integrals scale with a uniform mesh size $h$.

\subsection{Scaling Analysis in One Dimension}
\rev{We first consider a maximally smooth, untrimmed interior 1D basis function $N_i(\xi)$ with support $\mathcal{S}_i$ on a uniform knot vector. Its magnitude remains $\mathcal{O}(1)$ (i.e., $h^0$) under refinement, while its derivative scales as $\mathcal{O}(h^{-1})$ due to the shrinking support width.}

\paragraph{Scaling of the Mass Entry ($M_{ii} \sim h^1$)}
The property of partition of unity ensures that the sum of all basis functions reproduces the constant function. Integrating this identity over a single element of length $h$ yields:
\begin{equation}
\int_0^h \sum_j N_j(\xi) \, \mathrm{d}\xi
=
\int_0^h 1 \, \mathrm{d}\xi
=
h .
\end{equation}
\rev{For maximally smooth interior B-splines on a uniform knot vector, translational invariance implies that each basis function is a shifted version of the same reference spline. The $p+1$ contributions of the active basis functions over one element correspond, up to translation, to the $p+1$ element-wise contributions of a single interior basis function over its complete support $\mathcal{S}_i$. Consequently,}
\begin{equation}
M_{ii}
=
\int_{\mathcal{S}_i} N_i(\xi) \,\mathrm{d}\xi
=
\int_0^h \sum_j N_j(\xi) \,\mathrm{d}\xi
=
h
\sim h^1 .
\label{eq:app_mass_scaling_1d}
\end{equation}
\rev{The same result follows directly from the standard integral identity for a B-spline basis function \cite{LycheMorken2018},}
\begin{equation}
\int_{\mathcal{S}_i} N_i(\xi)\,\mathrm{d}\xi
=
\frac{\xi_{i+p+1}-\xi_i}{p+1}.
\end{equation}
\rev{For a maximally smooth interior B-spline on a uniform knot vector, $\xi_{i+p+1}-\xi_i=(p+1)h$, and hence}
\begin{equation}
\int_{\mathcal{S}_i} N_i(\xi)\,\mathrm{d}\xi
=
h,
\end{equation}
\rev{which confirms Eq.~\ref{eq:app_mass_scaling_1d}. This identity does not directly apply to the active portion of a trimmed basis function, for which the integration domain is reduced.}

\paragraph{Scaling of the Stiffness Entry ($K_{ii} \sim h^{-1}$)}
By the chain rule, the derivative scales inversely with the mesh size, $N_i' \sim h^{-1}$. Integrating the squared derivative over the support yields:
\begin{equation}
K_{ii}
=
\int_{\mathcal{S}_i}
\underbrace{[N_i'(\xi)]^2}_{\sim h^{-2}}
\underbrace{\mathrm{d}\xi}_{\sim h^1}
\sim h^{-1}.
\end{equation}
The ratio of these terms defines the Rayleigh quotient $Q_{ii}$. Substituting the scaling results, we obtain:
\begin{equation}
Q_{ii}
=
\frac{K_{ii}}{M_{ii}}
\sim
\frac{h^{-1}}{h^1}
=
h^{-2}
\implies
\omega_{i}^{\mathrm{Q}}
\sim
\sqrt{h^{-2}}
=
h^{-1}.
\end{equation}
This confirms that the local frequency is inversely proportional to the element length. We now show that this result is independent of the spatial dimension $d$.

\subsection{Generalization to Multiple Dimensions}
Multivariate B-splines are constructed as tensor products
$N_{\mathbf{i}}(\boldsymbol{\xi}) = \prod_{k=1}^d N_{i_k}(\xi_k)$.
Let $\mathcal{S}_i$ denote the $d$-dimensional support and $\mathcal{S}_{i_k}$ the 1D support in each direction.

\paragraph{Scaling of the Mass Entry ($M_{ii} \sim h^d$)}
Using the Theorem of Fubini, the $d$-dimensional mass integral decomposes into a product of 1D integrals. Since each univariate integral scales with $h^1$, the diagonal mass entry scales with the measure (area or volume) of the support:
\begin{equation}
M_{ii}
=
\int_{\mathcal{S}_i}
N_{\mathbf{i}}(\boldsymbol{\xi})
\,\mathrm{d}\boldsymbol{\xi}
=
\prod_{k=1}^d
\underbrace{
\int_{\mathcal{S}_{i_k}}
N_{i_k}(\xi_k)\,\mathrm{d}\xi_k
}_{\sim h^1}
\sim h^d .
\end{equation}

\paragraph{Scaling of the Stiffness Entry ($K_{ii} \sim h^{d-2}$)}
The stiffness matrix involves the squared norm of the gradient
\[
\|\nabla N_{\mathbf{i}}\|^2
=
\sum_{k=1}^d
\left(
\frac{\partial N_{\mathbf{i}}}{\partial \xi_k}
\right)^2 .
\]
A single partial derivative scales as:
\begin{equation}
\frac{\partial N_{\mathbf{i}}}{\partial \xi_k}
=
\underbrace{
\underbrace{
\frac{\partial N_{i_k}(\xi_k)}{\partial \xi_k}
}_{\sim h^{-1}}
\cdot
\underbrace{
\prod_{j \neq k}^d N_{i_j}(\xi_j)
}_{\sim h^{0}}
}_{\sim h^{-1}} .
\end{equation}
Integrating the squared partial derivative over the $d$-dimensional support introduces the volumetric measure $\mathrm{d}\boldsymbol{\xi}\sim h^d$:
\begin{equation}
\int_{\mathcal{S}_i}
\underbrace{
\left(
\frac{\partial N_{\mathbf{i}}}{\partial \xi_k}
\right)^2
}_{\sim h^{-2}}
\underbrace{
\mathrm{d}\boldsymbol{\xi}
}_{\sim h^d}
\sim h^{d-2}.
\end{equation}
The total diagonal stiffness entry $K_{ii}$ is the sum of $d$ such integrals, thus maintaining the scaling $\mathcal{O}(h^{d-2})$.

\subsection{Scaling of the Maximum Eigenfrequency}
The preceding DOF-wise analysis shows that the diagonal stiffness and row-sum-lumped mass contributions associated with an interior basis function scale as
\[
K_{ii}\sim h^{d-2},
\qquad
M_{ii}\sim h^d.
\]
Consequently, the corresponding diagonal Rayleigh quotient satisfies
\[
Q_{ii}
=
\frac{K_{ii}}{M_{ii}}
\sim h^{-2}.
\]
Equivalently, for a fixed polynomial degree and a fixed type of interior basis function, one may write
\[
Q_{ii}=C_Q h^{-2},
\]
where the constant $C_Q>0$ is independent of the mesh size $h$.
The maximum diagonal Rayleigh quotient provides a lower bound for the largest discrete eigenvalue,
\[
\max_i Q_{ii}\le\lambda_{\max}.
\]
\rev{Since the diagonal Rayleigh quotient of an interior basis function scales as $Q_{ii}=C_Qh^{-2}$, and the maximum over all DOFs cannot be smaller than this value, it follows that}
\[
C_Q h^{-2}
\le
\max_i Q_{ii}
\le
\lambda_{\max}.
\]
\rev{Hence, the DOF-wise Rayleigh analysis shows that the largest eigenvalue is bounded from below by a quantity proportional to $h^{-2}$.}

\rev{Let $v_h=\sum_i q_i N_i$ denote the discrete function associated with the coefficient vector $\mathbf q$. To obtain the corresponding upper scaling, we invoke the standard inverse inequality for fixed polynomial degree on shape-regular discretizations \cite{Gallistl2017}. This is a classical estimate and is used here only to complete the scaling argument for $\lambda_{\max}$.}
\[
\|\nabla v_h\|_{L^2}
\le
C_{\mathrm{inv}}h^{-1}
\|v_h\|_{L^2},
\]
\rev{where $C_{\mathrm{inv}}$ is independent of $h$. Squaring this relation gives}
\[
\|\nabla v_h\|_{L^2}^2
\le
C_{\mathrm{inv}}^2 h^{-2}
\|v_h\|_{L^2}^2.
\]
\rev{For the discrete problem considered here, the numerator of the Rayleigh quotient represents the stiffness contribution and is proportional to $\|\nabla v_h\|_{L^2}^2$, while the denominator represents the mass contribution. For fixed polynomial degree on a uniform mesh, the norm induced by the row-sum-lumped mass matrix is equivalent to the corresponding $L^2$ norm. Therefore, the Rayleigh quotient of any admissible discrete function satisfies}
\[
\frac{\mathbf q^T\mathbf K\mathbf q}
     {\mathbf q^T\mathbf M\mathbf q}
\le
C\,h^{-2},
\]
\rev{with a constant $C>0$ independent of $h$. Since the largest eigenvalue is the maximum of this quotient over all non-zero vectors $\mathbf q$, it follows that}
\[
\lambda_{\max}
=
\max_{\mathbf q\neq\mathbf 0}
\frac{\mathbf q^T\mathbf K\mathbf q}
     {\mathbf q^T\mathbf M\mathbf q}
\le
C\,h^{-2}.
\]
\rev{The lower and upper estimates therefore have the same dependence on the mesh size,}
\[
C_Q h^{-2}
\le
\lambda_{\max}
\le
C h^{-2}.
\]
\rev{Since both constants are independent of $h$, the maximum eigenvalue exhibits the familiar asymptotic behavior}
\begin{equation}
\lambda_{\max}(h)\sim h^{-2}.
\end{equation}
\rev{Since $\lambda_{\max}=\omega_{\max}^2$, the corresponding maximum eigenfrequency scales as}
\begin{equation}
\omega_{\max}(h)\sim h^{-1}.
\end{equation}
\rev{Thus, under uniform bisection, $h\rightarrow h/2$,}
\[
\omega_{\max}(h/2)
\sim
2\,\omega_{\max}(h),
\]
\rev{and the maximum eigenfrequency is asymptotically doubled.}

\paragraph{Remark (limitations)}
\rev{The above scaling argument is a standard result for maximally smooth, untrimmed interior basis functions on uniform or shape-regular meshes with fixed polynomial degree. It is included here only to make the mesh-size dependence used in the main text explicit. The argument does not directly extend to arbitrarily trimmed basis functions, since trimming changes the active support and may modify stiffness and row-sum-lumped mass contributions non-proportionally. The resulting frequency bounds can therefore depend on the particular trimming configuration, as investigated in Sec.~\ref{subsub:Rayleigh}.}

\clearpage

\section*{CRediT authorship contribution statement}
\noindent 
\textbf{Christoph Hollweck:} Software, Validation, Investigation, Conceptualization, Writing - Original Draft, Visualization.\\
\textbf{Lukas Leidinger:}  Software, Validation, Investigation, Conceptualization, Writing - review and editing.\\
\textbf{Stefan Hartmann:} Software, Validation, Investigation.\\ 
\textbf{Marcus Wagner:} Review, Supervision.\\
\textbf{Roland Wüchner:} Review, Supervision.

\section*{Declaration of competing interest}
\noindent
The authors declare no competing financial interests. The results of this study contribute to ongoing developments within the commercial solver LS-DYNA.

\section*{Data availability}
\noindent 
Data will be made available on request.

\section*{Declaration of generative AI and AI-assisted technologies in the writing process}
\noindent During the preparation of this work, the authors used large language models (LLMs) solely to improve readability and language. The authors reviewed and edited the content as needed and take full responsibility for the content of the publication.

\section*{Acknowledgements}
\noindent C. Hollweck is supported by a PhD scholarship from the Studienstiftung des deutschen Volkes and by OTH Regensburg. The authors gratefully acknowledge this support.


\clearpage
\bibliographystyle{elsarticle-num} 
\bibliography{mylib}





\end{document}